\documentclass[11pt]{article}
\usepackage[margin=1in]{geometry}
\usepackage{amsmath,amssymb,amsthm,mathtools,mathrsfs,booktabs,longtable,microtype,needspace}
\usepackage[hidelinks]{hyperref}
\hypersetup{pdftitle={Circular flows and bisections of random regular graphs},pdfauthor={Amir Jafari},pdfsubject={Circular flows and random regular graphs},pdfkeywords={circular flow, random regular graph, prescribed outdegree, maximum bisection, entropy method}}
\newtheorem{theorem}{Theorem}[section]
\newtheorem{lemma}[theorem]{Lemma}
\newtheorem{proposition}[theorem]{Proposition}
\newtheorem{corollary}[theorem]{Corollary}
\newtheorem{conjecture}[theorem]{Conjecture}
\newcommand{\E}{\mathbb E}
\newcommand{\cut}{\partial}
\newcommand{\bis}{\operatorname{bis}}
\numberwithin{equation}{section}
\title{Circular flows and bisections of random regular graphs}
\author{Amir Jafari\\[2pt]\small New Uzbekistan University}\date{}
\begin{document}
\maketitle
\begin{abstract}
For every fixed odd $d\ge5$, we prove that the circular flow number of a uniformly random simple $d$-regular graph is asymptotically almost surely less than $2+4/(d-1)-\varepsilon_d$ for an explicit constant $\varepsilon_d>0$. At degree five the upper bound is $3-1/24006$. At degree seven it is $21/8$, answering a question of Li and Li. For every fixed $p\ge1$, we also obtain a strongly connected Jaeger orientation at degree $4p+1$, with outdegrees in $\{p,3p+1\}$. At large degrees we count dense bisections together with sets that violate their cut inequalities. This gives flow intervals of length $O(d^{-2})$ above the maximum-bisection lower bound of Li and Li. The finite-degree estimates, together with the large-degree transfer, also prove the conjecture of Greenhill, Isaev and Lewis that the threshold for orientations with outdegrees in $\{p,d-p\}$ exceeds the corresponding bisection threshold by at most one.
\end{abstract}

\section{Introduction}\label{sec:introduction}

Li and Li \cite[Theorem~1.4]{LL} determined the leading asymptotic behaviour of the circular flow number of random regular graphs of large odd degree. We prove a fixed gap below $2+4/(d-1)$ at every odd degree $d\ge5$ and sharpen the large-degree estimate in terms of the maximum bisection of the sampled graph.

At degree $d=4p+1$, the value $2+4/(d-1)=2+1/p$ comes from Jaeger's circular flow conjecture. For a positive integer $p$, Jaeger conjectured that every $4p$-edge-connected graph has an orientation whose outdegree and indegree at each vertex are congruent modulo $2p+1$, equivalently a circular flow number at most $2+1/p$ \cite[Section~1]{LTWZ}. For $(4p+1)$-regular graphs these are the \emph{Jaeger orientations}, with outdegrees $p$ and $3p+1$. Han, Li, Wu and Zhang \cite[Theorem~1.3]{HLWZ} constructed counterexamples to the general conjecture for every $p\ge3$. In the random regular model, Pra\l{}at and Wormald \cite[Theorem~1.2]{PW} proved the case $p=1$, Delcourt, Huq and Pra\l{}at \cite[Theorem~4]{DHP} proved $p=2$, and Li, Li and Liu \cite[Theorem~1.3]{LLL} established all cases $p\le4$. Alon and Pra\l{}at \cite[Theorem~1.6]{AP} proved existence for all sufficiently large fixed $p$.

Write $G_{n,d}$ for a uniformly random simple $d$-regular graph on $[n]=\{1,\ldots,n\}$. Throughout, $d\ge3$ is fixed. Except in Section~\ref{sec:questions}, $n$ tends to infinity through even integers. A property holds asymptotically almost surely, abbreviated a.a.s., if its probability tends to one. All logarithms are natural. An unsubscripted $O(\cdot)$ has an absolute implied constant unless other fixed parameters are specified; $O_d(\cdot)$ allows dependence on $d$. The notation $o_n(1)$ means convergence to zero with the other parameters fixed. Large-degree limits, denoted $o_d(1)$ or $o_k(1)$, are taken after the fixed-degree probabilistic statement.

For a graph $G$, put $n=|V(G)|$. For $S\subseteq V(G)$, let $\cut(S)$ be the edges with exactly one endpoint in $S$. A nontrivial set is nonempty and proper. A bisection partitions $V(G)$ into two equal parts. Write
\[
 \bis(G)=\max\{\text{number of edges between the parts of a bisection}\},
 \qquad c(G)=\bis(G)/n.
\]
Thus $c(G)$ is measured per vertex. A balanced signing is a map $\sigma:V(G)\to\{-1,1\}$, with $\sigma_v=\sigma(v)$, such that $\sum_v\sigma_v=0$. Put
\[
 \sigma(S)=\sum_{v\in S}\sigma_v,\qquad
 Q_G(\sigma)=|\cut(\sigma^{-1}(1))|,\qquad c_\sigma=Q_G(\sigma)/n.
\]
For disjoint $U,W$, write $E(U,W)$ for the edges joining them.

An orientation is strongly connected if every ordered pair of vertices is joined by a directed path. Its divergence at $v$ is its outdegree minus its indegree. A real circulation assigns edge values with weighted divergence zero at every vertex. A circular $r$-flow is an orientation with a circulation taking values in $[1,r-1]$. The circular flow number $\phi(G)$ is the infimum of such $r\ge2$, with value $+\infty$ if none exists.

\begin{theorem}\label{main:jaeger}
For every fixed integer $p\ge1$, there is a constant $\eta_p>0$ such that, a.a.s., $G_{n,4p+1}$ has a strongly connected Jaeger orientation and
\[
 \phi(G_{n,4p+1})<2+1/p-\eta_p.
\]
One may take $\eta_1=\varepsilon_5=1/24006$.
\end{theorem}

In general one may take $\eta_p=\varepsilon_{4p+1}$ from Theorem~\ref{main:odd-flow}.

The orientation assertion proves \cite[Conjecture~1.2]{GIL} and answers \cite[Problem~1.5]{LLL}. The cut estimate of Alon and Pra\l{}at \cite[Theorem~1.5]{AP} gives strong connectivity and a positive flow margin at sufficiently large degree.

The circular-flow conclusion extends to all odd degrees at least five.
\begin{theorem}\label{main:odd-flow}
For every fixed odd $d\ge5$, there exists $\varepsilon_d>0$, independent of $n$, such that
\[
 \phi(G_{n,d})<2+\frac4{d-1}-\varepsilon_d
 \qquad\text{a.a.s.}
\]
The following explicit bounds hold a.a.s.
\[
\begin{array}{c|c|c}
 d&\text{strict upper bound for }\phi(G_{n,d})&\varepsilon_d\\\hline
 5&3-1/24006&1/24006\\
 7&21/8&1/24\\
 9&202/81&1/162\\
 11&1202/501&2/2505\\
 13&142/61&1/183\\
 15&1602/701&2/4907\\
 17&92/41&1/164\\
 19&202/91&2/819
\end{array}
\]
For every fixed odd $d\ge21$, one also has
\[
 \phi(G_{n,d})<2+\frac4{d-6+(4/\pi)\sqrt d}
 <2+\frac4{d-1}.
\]
For these degrees one may take $\varepsilon_d$ to be the difference between the last two bounds.
\end{theorem}

Li and Li \cite[Proposition~5.1]{LL} proved $\phi(G_{n,7})<3$ a.a.s. and asked for a fixed upper bound below $3$ \cite[Section~5, p.~10]{LL}. The bound $21/8$ answers this question and lies $1/24$ below $2+4/(d-1)=8/3$. Their large-degree asymptotic \cite[Theorem~1.4]{LL} also implies a fixed gap for all sufficiently large odd degrees.

\subsection{Bisections and the orientation threshold}

For an integer $1\le p\le\lfloor d/2\rfloor$, a $p$-orientation has every outdegree in $\{p,d-p\}$. When $p<d/2$, the two vertex types have divergences $\pm(d-2p)$ and each contains $n/2$ vertices. Summing divergence over one type gives the necessary condition $c(G)\ge d/2-p$. Following \cite[Section~1.1]{GIL}, define
\begin{align*}
 p^*(d)&=\min\{p:1\le p\le\lfloor d/2\rfloor,\quad
              G_{n,d}\text{ a.a.s. has a $p$-orientation}\},\\
 p_{\rm bis}(d)&=\min\{p:1\le p<d/2,\quad
              G_{n,d}\text{ a.a.s. has }c(G_{n,d})\ge d/2-p\},
\end{align*}
where the minima are over integers and an empty minimum is $+\infty$. For a sampled graph, let $p_{\rm or}(G)$ be its least orientation parameter, with the same convention.

\begin{theorem}\label{main:gap}
For every integer $d\ge3$,
\[
 p_{\rm bis}(d)\le p^*(d)\le p_{\rm bis}(d)+1.
\]
For each fixed $d\ge21$, a.a.s.
\begin{equation}\label{eq:pointwiseintro}
 0\le p_{\rm or}(G_{n,d})-\lceil d/2-c(G_{n,d})\rceil\le1.
\end{equation}
\end{theorem}

The first assertion proves \cite[Conjecture~1.5]{GIL}. The case $d=8$ follows by adding a matching to the degree-seven construction, which uses the preprint count of Greenhill, Isaev and Lewis \cite{GIL}. The other degrees in this theorem use published inputs and the estimates proved here.

\subsection{Circular-flow intervals}

Li and Li \cite[Lemma~2.4]{LL} proved that every loopless odd-regular graph of degree at least three satisfies
\begin{equation}\label{eq:lili-lower}
 \phi(G)\ge\frac{4c(G)}{2c(G)-1}=2+\frac2{2c(G)-1}.
\end{equation}
Lemma~\ref{s1:lem:criterion} gives a common cut-ratio formulation and includes a short proof of this lower bound. Our upper bounds measure the loss from this obstruction in terms of the sampled graph's actual maximum bisection.

\begin{theorem}\label{main:circular}
For every fixed odd $d\ge21$, a.a.s., with $c=c(G_{n,d})$,
\begin{equation}\label{eq:loss2intro}
 2+\frac2{2c-1}\le\phi(G_{n,d})<2+\frac2{2c-3}.
\end{equation}
For every fixed odd $d\ge2501$, the upper bound improves to
\begin{equation}\label{eq:loss1intro}
 \phi(G_{n,d})<2+\frac2{2c-2}.
\end{equation}
For every sufficiently large fixed odd $d$, it improves further to
\begin{equation}\label{eq:halfwindowintro}
 2+\frac2{2c-1}\le\phi(G_{n,d})<2+\frac2{2c-3/2}.
\end{equation}
When $d=2k+1$, the interval in \eqref{eq:halfwindowintro} has length $(1+o_k(1))k^{-2}$. Both endpoints use the maximum bisection density of the sampled graph.
\end{theorem}

Li and Li \cite[Theorem~3.7]{LL} transfer a large bisection to cut inequalities. Our transfer holds uniformly over dense signings, with explicit degree thresholds and constant losses. Section~\ref{s2:sec:precision} also refines the remainder in their large-degree asymptotic, using the bisection estimate proved in Appendix~\ref{app:dms}.

For degree $4p+1$, the strict bound $\phi(G)<2+1/p$ is equivalent to the existence of a strongly connected orientation with divergence divisible by $2p+1$ \cite[Theorem~1.1]{LTWZ}. Jafari \cite[Corollary~1.3 and Proposition~4.5]{GB} proves this conclusion for every $(4p+2)$-edge-connected graph. A cut isolating a vertex of a $(4p+1)$-regular graph has size $4p+1$, and strictness can fail in this degree, as $\phi(K_{4p+2})=2+1/p$ \cite[Section~5, p.~176]{LTWZ}.

\subsection{The argument}

Our main estimate transfers bisection density to cut inequalities. For each fixed $d\ge21$, a.a.s. every balanced signing of $G_{n,d}$ with $c_\sigma\ge d/4+(19/60)\sqrt d$ satisfies
\[
 |\cut(S)|>(2c_\sigma-2)|\sigma(S)|
 \qquad(\varnothing\ne S\ne V(G_{n,d})).
\]
Lyons's bisection theorem \cite[Theorem~4.1, equation~(4.2)]{LY} supplies such a signing, and the cut criterion gives the bound of Theorem~\ref{main:odd-flow} for $d\ge21$, together with Jaeger orientations for $p\ge5$.

A maximum bisection is selected from the sampled graph, so an estimate for one fixed signing does not control its cuts. We therefore count a bisection together with a set violating its cut inequality. Requiring density at least $d/4+\beta\sqrt d$ reduces the leading first-moment exponent for signings from $\log2$ to $\log2-4\beta^2$ as $d$ grows. For the localized joint count in Section~\ref{s2:sec:precision}, a loss $\ell$ in the cut parameter contributes a further term $-\tfrac\ell2\log T_0$ to the limiting exponent, where $T_0$ is the square of the tribonacci constant.

The four intersections of the bisection parts with the set and its complement record both cuts. Comparison weights cancel the small-class logarithms. A quadratic entropy bound handles sets close to a bisection side, and convexity reduces the remaining region to finitely many interval endpoints in Appendix~\ref{s3:app:tables}. Harangi \cite[Lemma~2.1, Sections~3--4]{H} also uses first-moment counts with intersections of prescribed classes. Here the joint count includes the bisection and its crossing-density constraint, giving uniformity over sufficiently dense bisections.

The smaller degrees require additional constructions. At degrees nine, thirteen and seventeen, we transfer a directed cut bound from a model weighted by its number of Jaeger orientations. Adding a perfect matching to a six-regular factor gives degree seven, and adding a Hamilton cycle to a nine-regular factor gives degree eleven. The dense-bisection count also covers degrees fifteen and nineteen. At degree five, we construct bounded monochromatic trees and prove expansion of their component graph.

Section~\ref{s1:sec:jaeger} gives the cut-ratio criterion and applies it to the flow bounds. Section~\ref{sec:counting} proves the dense-bisection transfer, Section~\ref{s3:sec:gap} deduces the one-unit gap, and Section~\ref{s2:sec:precision} sharpens the transfer at large degrees. Section~\ref{sec:degreefive} proves the separate degree-five bound.
\section{Strict circular flows from dense bisections}\label{s1:sec:jaeger}

For a balanced signing $\sigma$ of a graph of even order, put
\[
 q_G(\sigma)=\min_{\sigma(S)\ne0}\frac{|\cut(S)|}{|\sigma(S)|},
 \qquad q^*(G)=\max_\sigma q_G(\sigma).
\]
The minima and maximum are over finite nonempty sets. Since the positive
side of $\sigma$ is a bisection, $q_G(\sigma)\le2c_\sigma$ and
$q^*(G)\le2c(G)$. We call $q_G(\sigma)$ the cut ratio of the signing and $q^*(G)$ the cut ratio of the graph. The following criterion combines the prescribed-degree and circulation conditions. The balanced-signing formulation appears in Li and Li \cite[Lemma~2.3]{LL}, and we allow any real flow parameter.

\begin{lemma}[The cut-ratio criterion]\label{s1:lem:criterion}\label{prop:cut-ratio}
Let $G$ be a $d$-regular graph of even order, with $d\ge3$, and let $\sigma$ be a balanced signing.
\begin{enumerate}
\item Let $q>0$ be an integer of the parity of $d$. An orientation with
      divergence $q\sigma$ exists if and only if $q\le q_G(\sigma)$.
      If $G$ is connected, such an orientation is strongly connected if
      and only if $q<q_G(\sigma)$. Consequently
\[
 p_{\rm or}(G)=\max\left\{1,\left\lceil\frac{d-q^*(G)}2\right\rceil\right\}.
\]
\item If $d$ is odd and $q_G(\sigma)>1$, then
\[
 \phi(G)\le 2+\frac2{q_G(\sigma)-1}.
\]
      In particular, $q_G(\sigma)>q>1$ gives
      $\phi(G)<2+2/(q-1)$.
\item Suppose $d$ is odd, and put $r_d=2d/(d-3)$ for $d\ge5$ and
      $r_3=+\infty$. Then
\[
 \phi(G)<r_d\quad\Longleftrightarrow\quad q^*(G)>d/3.
\]
      Under these equivalent conditions,
\[
 \phi(G)=2+\frac2{q^*(G)-1}.
\]
      If $G$ is loopless, the lower bound of Li and Li follows for every odd $d\ge3$:
\[
 \phi(G)\ge2+\frac2{2c(G)-1}.
\]
\end{enumerate}
\end{lemma}

\begin{proof}
Hakimi's prescribed-degree criterion, in the form
\cite[Theorem~3]{DGKNS}, gives the first equivalence. The required
outdegrees are $(d+q\sigma_v)/2$. Writing $e(S)$ for the number of edges with both endpoints in $S$, the set condition is
$\sum_{v\in S}(d+q\sigma_v)/2-e(S)=\{|\cut(S)|+q\sigma(S)\}/2\ge0$. A nontrivial directed cut is
excluded precisely when $|\cut(S)|>q|\sigma(S)|$; for connected $G$,
sets with $\sigma(S)=0$ already satisfy strictness. The formula for
$p_{\rm or}$ follows by taking $q=d-2p$. At even $d$, the remaining
case $p=d/2$ is Eulerian. At odd $d$, $q^*(G)\ge1$: add any perfect
matching on the vertex set, orient the resulting even-degree graph
Eulerian, and delete the matching.

For the second assertion, take an orientation $D$ of divergence
$\sigma$, supplied by the first assertion, and put
$r=2+2/(q_G(\sigma)-1)$. Writing $d_D^+(S)$ and $d_D^-(S)$ for the numbers of edges leaving and entering $S$, Hoffman's circulation criterion
\cite[Theorem~1]{DGKNS} gives a circulation in $[1,r-1]$, because
\[
 (r-1)d_D^-(S)-d_D^+(S)
 =\frac{r-2}{2}|\cut(S)|-\frac r2\sigma(S)\ge0.
\]

Conversely, let a circular $r$-flow, with $r>2$, have orientation divergence $b$, and put $b(S)=\sum_{v\in S}b(v)$.
Subtracting $r/2$ from every edge value and summing over a set gives
\[
 \frac r{r-2}|b(S)|\le|\cut(S)|.
\]
If $r<r_d$, then $r/(r-2)>d/3$. The singleton inequalities give
$|b(v)|<3$, and odd degree forces $b(v)=\pm1$. Thus $b$ is a balanced
signing and $q^*(G)\ge r/(r-2)$. Together with the second assertion,
this proves the equivalence. Passing to a sequence $r\downarrow\phi(G)$
gives the stated identity.
Finally, if $G$ is loopless, averaging over bisections gives $2c(G)\ge d/2>d/3$.
If $\phi(G)<2+2/(2c(G)-1)$, the equality just proved and
$q^*(G)\le2c(G)$ give a contradiction.
\end{proof}

\subsection{Flow margins from orientation counts}

For $p=2,3,4$, the normalized orientation count of Li, Li and Liu converges in distribution to a positive random variable. This allows us to transfer a cut inequality from a marked orientation to almost every Jaeger orientation of a typical graph.

In the marked model, indegrees and outdegrees are bounded and at least two. The diameter theorem of Caputo and Quattropani \cite[Theorem~1.1]{CQ} gives strong connectivity under these hypotheses. The directed cut estimate in Appendix~\ref{app:marked-expansion} supplies the quantitative margin for our flow bound.

\begin{lemma}[Transfer from the marked model]\label{lem:marked-transfer}
For each $n$, let $X_n$ and $B_n$ be random variables with $0\le B_n\le X_n$ and $\mu_n=\E X_n>0$. Suppose that
\[
 X_n/\mu_n\xrightarrow{\mathrm d}W,\qquad
 \mathbb P(W>0)=1,\qquad \E B_n/\mu_n\longrightarrow0,
\]
where $\xrightarrow{\mathrm d}$ denotes convergence in distribution. Then $X_n>0$ a.a.s., and $B_n/X_n\to0$ in probability, the ratio being taken as zero when $X_n=0$. In particular, $X_n>B_n$ a.a.s.
\end{lemma}

\begin{proof}
For $\rho,\eta>0$, Markov's inequality gives
\[
 \mathbb P(B_n>\rho X_n)
 \le \mathbb P(X_n\le\eta\mu_n)
       +\frac{\E B_n}{\rho\eta\mu_n}.
\]
The limit superior is at most $\mathbb P(W\le\eta)$. Letting $\eta\downarrow0$ proves the ratio assertion. The same lower-tail bound gives $\mathbb P(X_n=0)\to0$. Taking $\rho=1/2$ now gives $X_n>B_n$ a.a.s.
\end{proof}

In our applications $X_n$ counts the marked structures of a random pairing, possibly together with an independent auxiliary structure, and $B_n$ counts those that fail a required property. Then $\E B_n/\E X_n$ is the failure probability under the marked law, which weights each outcome of the base law by $X_n/\E X_n$ and then chooses one of its marks uniformly.

\begin{lemma}[Directed cuts and flow margins]\label{lem:directed-margin}
Let $G$ be a $d$-regular graph on $n$ vertices, and let $D$ be an orientation of $G$ with divergence $q\sigma$, where $q>1$ and $\sigma$ is a balanced signing. Suppose that, for some $h>0$,
\[
 \min\{d_D^+(S),d_D^-(S)\}
 >h\min\{|S|,n-|S|\}
 \qquad(\varnothing\ne S\ne V(G)).
\]
Then $D$ is strongly connected, and $|\cut(S)|>(q+2h)|\sigma(S)|$ for every nontrivial $S$. If $d$ is odd, then
\[
 \phi(G)<2+\frac2{q+2h-1}
 =2+\frac2{q-1}-\varepsilon(q,h),\qquad
 \varepsilon(q,h)=\frac{4h}{(q-1)(q+2h-1)}>0.
\]
\end{lemma}
\begin{proof}
The hypothesis excludes directed cuts, so $D$ is strongly connected. Balance gives $|\sigma(S)|\le\min\{|S|,n-|S|\}$. Summing divergence over $S$ therefore yields
\begin{align*}
 |\cut(S)|
 &=q|\sigma(S)|+2\min\{d_D^+(S),d_D^-(S)\}\\
 &>q|\sigma(S)|+2h\min\{|S|,n-|S|\}
 \ge(q+2h)|\sigma(S)|.
\end{align*}
Thus $q_G(\sigma)>q+2h$, and Lemma~\ref{s1:lem:criterion}(2) gives the flow bound.
\end{proof}

\begin{proposition}\label{prop:small-strong}
Put $h_2=1/40$, $h_3=1/20$ and $h_4=1/10$. For each $p\in\{2,3,4\}$, a.a.s. $G_{n,4p+1}$ has a Jaeger orientation $D$ such that
\begin{equation}\label{eq:uniform-directed-cuts}
 \min\{d_D^+(S),d_D^-(S)\}>h_p\min\{|S|,n-|S|\}
 \qquad(\varnothing\ne S\ne V(G)).
\end{equation}
The proportion of its Jaeger orientations satisfying this inequality tends to one in probability, taking the proportion to be zero when there is no Jaeger orientation. Consequently, a.a.s.
\begin{equation}\label{eq:small-fixed-margin}
 \phi(G_{n,4p+1})<2+\frac1{p+h_p}
 =2+\frac1p-\frac{h_p}{p(p+h_p)}.
\end{equation}
\end{proposition}
\begin{proof}
Fix $p\in\{2,3,4\}$ and put $d=4p+1$. In the pairing model, each vertex has $d$ labelled points, and a uniformly chosen perfect matching of all $dn$ points gives a $d$-regular multigraph. Let $X_n$ count its Jaeger orientations and let $B_n$ count those failing~\eqref{eq:uniform-directed-cuts}. Orientations here distinguish the two endpoint points of every pair, including loops.

Li, Li and Liu \cite[Lemmas~2.2--2.3, Theorem~2.4 and Appendix~A]{LLL} verify the small-subgraph-conditioning hypotheses for $X_n$. Janson's distributional conclusion, in the form \cite[Theorem~4.3]{W}, therefore gives
\begin{equation}\label{eq:orientation-ssc}
 \frac{X_n}{\E X_n}\xrightarrow{\mathrm d}W_p
 =\prod_{j\ge1}(1+\delta_j)^{Z_j}e^{-\lambda_j\delta_j},
 \qquad
 \lambda_j=\frac{(4p)^j}{2j},\quad
 \delta_j=\left(-\frac p{4p+1}\right)^j,
\end{equation}
where the $Z_j$ are independent Poisson random variables with means $\lambda_j$. Every $\delta_j>-1$, and
\[
 \sum_{j\ge1}\lambda_j\delta_j^2
 =-\frac12\log\left(1-\frac{4p^3}{(4p+1)^2}\right)<\infty.
\]
These conditions give $W_p>0$ almost surely; equivalently, the required lower-tail conclusion is \cite[Corollary~4.2]{W}.

Now choose uniformly a pair consisting of a pairing and one of its Jaeger orientations. To generate this law, proceed as follows. Choose the $n/2$ vertices with outdegree $p$, choose the outgoing points at every vertex, and match all outgoing points bijectively to all incoming points, uniformly. There are
\[
 \binom n{n/2}\binom d p^{\!n}(dn/2)!
\]
equally likely outcomes, one for each oriented pairing. Conditional on the vertex types and point choices, this is the directed configuration model with indegree and outdegree pairs $(p,3p+1)$ and $(3p+1,p)$. Lemma~\ref{lem:marked-expansion}, proved in Appendix~\ref{app:marked-expansion}, shows that it satisfies~\eqref{eq:uniform-directed-cuts} a.a.s.

Under this uniform law on oriented pairings, the probability of failure is exactly
\[
 \frac{\sum_P B_n(P)}{\sum_P X_n(P)}
 =\frac{\E B_n}{\E X_n}=o(1),
\]
where the sums range over all pairings $P$. Every Jaeger orientation has divergence $(2p+1)\sigma$ for a balanced signing $\sigma$. Lemma~\ref{lem:marked-transfer} gives the assertions about the number and proportion of orientations satisfying~\eqref{eq:uniform-directed-cuts} in the uniform pairing model. Lemma~\ref{lem:directed-margin}, with $q=2p+1$ and $h=h_p$, then gives~\eqref{eq:small-fixed-margin}. Conditional on simplicity the graph is uniform in $G_{n,d}$, and the conditioning event has a positive limiting probability \cite[Section~2.2, equation~(3)]{W}. Thus all assertions pass to simple graphs.

\end{proof}

The directed cut estimate requires $p\ge2$. The degree-five case uses the recolouring argument in Section~\ref{sec:degreefive} and Appendix~\ref{app:p1strict}.

\subsection{Degree seven}

\begin{proposition}\label{prop:degree7}
Asymptotically almost surely, $q^*(G_{n,7})>21/5$ and $\phi(G_{n,7})<21/8=8/3-1/24$.
\end{proposition}
\begin{proof}
Start with a six-regular pairing $P$ marked by an orientation with all outdegrees one or five. There are $n/2$ of each type. Let $\sigma$ be positive on the vertices of outdegree five. Adjoin an independent uniform perfect matching $M_0$ of the vertex set. Under the uniform law on marked pairings, conditional on the types and point directions, the outgoing points are matched uniformly to the incoming points. Lemma~\ref{lem:degree-seven-marked} shows that with probability $1-o(1)$ this signing satisfies $|\cut_{P\cup M_0}(S)|>(21/5)|\sigma(S)|$ on every nontrivial set.

Now take $P$ uniform and independent of $M_0$. Let $X$ count all the marks of $P$ and $B$ count those failing the cut bound after adding $M_0$. The marked estimate is $\E B/\E X=o(1)$. Greenhill, Isaev and Lewis \cite[Lemmas~3.3--3.5 and~4.7]{GIL} verify the small-subgraph-conditioning hypotheses for $(d,p)=(6,1)$, with cycle parameters $\lambda_j=5^j/(2j)$ and $\delta_j=(-1/3)^j$. Since $\sum_j\lambda_j\delta_j^2<\infty$ and every $\delta_j>-1$, Janson's theorem \cite[Theorem~4.3]{W} gives $X/\E X\xrightarrow{\mathrm d}W$ with $W>0$ almost surely. Lemma~\ref{lem:marked-transfer} gives $X>B$ a.a.s.

The probability of a simple union is bounded away from zero, and conditional on simplicity its law is that of the union of independent uniform $G_{n,6}$ and a uniform perfect matching, conditioned to be simple (the graph-restricted sum $\mathcal G_{n,6}\oplus\mathcal G_{n,1}$ of \cite[Theorem~4.15(ii)]{W}). Janson's contiguity theorem \cite[Theorem~4.15(ii)]{W} transfers the existence of the signing to $G_{n,7}$. Now $q_G(\sigma)>21/5$, and Lemma~\ref{s1:lem:criterion}(2) gives $\phi(G)<2+2/(21/5-1)=21/8$.
\end{proof}

\subsection{Degree eleven}

\begin{proposition}\label{prop:degree11}
Asymptotically almost surely, $\phi(G_{n,11})<1202/501=12/5-2/2505$.
\end{proposition}
\begin{proof}
Put $(d,q,\beta)=(9,5,101/200)$. Let $P$ be a uniform degree-$d$ pairing and let $C$ be an independent uniform Hamilton cycle on the same vertices with a uniformly chosen cyclic direction. We show that a.a.s. $P$ has a Jaeger orientation whose extension to $P\cup C$ satisfies
\begin{equation}\label{eq:eleven-cut}
 |\cut(S)|>(q+2\beta)|\sigma(S)|
 \qquad(\varnothing\ne S\ne V(P)),
\end{equation}
where the original orientation has divergence $q\sigma$. The cycle adds zero divergence, and Lemma~\ref{s1:lem:criterion}(2) then gives the asserted bound.

Lemma~\ref{lem:eleven-marked} in Appendix~\ref{app:eleven-entropy} shows that, if a pair consisting of a degree-$d$ pairing and one of its Jaeger orientations is chosen uniformly and the independent cycle is adjoined, then \eqref{eq:eleven-cut} holds a.a.s.

To transfer this estimate, let $X_n(P)$ count the Jaeger orientations of $P$, and let $B_n(P,C)$ count those whose extension fails~\eqref{eq:eleven-cut}. The marked estimate gives $\E B_n(P,C)/\E X_n(P)=o(1)$. By~\eqref{eq:orientation-ssc} at degree nine,
\[
 X_n/\E X_n\xrightarrow{\mathrm d}W,
 \qquad \Pr(W>0)=1.
\]
Adjoining the independent cycle preserves this lower-tail bound. Lemma~\ref{lem:marked-transfer} therefore supplies an orientation satisfying~\eqref{eq:eleven-cut} a.a.s.

Finally condition on the union being simple, an event of probability bounded away from zero. Its conditional law is the union of an independent $G_{n,d}$ and uniform Hamilton cycle, conditioned on a simple union. This law and $G_{n,d+2}$ are mutually contiguous, meaning that they have the same a.a.s. events, by Janson's degree-augmentation theorem \cite{JAN}, in the form \cite[Theorem~4.15(i) and proof of Lemma~4.14]{W}. Hence the cut inequalities and flow bounds pass to the uniform simple graphs.
\end{proof}

\subsection{Dense bisections}

We use the following bisection bound of Lyons.
\begin{lemma}[Lyons]\label{s1:lem:bisectionstrong}\label{s1:lem:bisection}
For every fixed $d\ge3$, a.a.s.
\[
 c(G_{n,d})>\frac d{2\pi}\arccos\!\left(-\frac{2\sqrt{d-1}}d\right)
 >\frac d4+\frac{\sqrt d}{\pi}
 >\frac d4+\frac{19}{60}\sqrt d.
\]
\end{lemma}
\begin{proof}
Apply \cite[Theorem~4.1, equation~(4.2)]{LY} to an independent sequence $G_{n,d}$. As shown there, the proportion of vertices with a tree neighbourhood of each fixed radius tends to one almost surely. The resulting strict lower bound on the limit inferior gives the first inequality a.a.s. Equation~(4.2) also gives the second inequality. The last follows from $\pi<22/7<60/19$.
\end{proof}

Let $G$ be loopless and have odd degree $d\ge3$. Put $c=c(G)$ and fix $\ell>0$ with $2c-\ell>1$. Suppose a maximum-bisection signing satisfies the transfer inequality
\begin{equation}\label{s2:eq:transferhyp}
 |\cut(S)|>(2c-\ell)|\sigma(S)|
 \qquad(\varnothing\ne S\ne V(G)).
\end{equation}
Then $2c-\ell<q_G(\sigma)\le q^*(G)\le2c$. Lemma~\ref{s1:lem:criterion} therefore gives
\begin{equation}\label{s2:eq:window}
 2+\frac2{2c-1}\le\phi(G)<2+\frac2{2c-\ell-1}.
\end{equation}
The interval has length
\begin{equation}\label{s2:eq:width}
 \frac2{2c-\ell-1}-\frac2{2c-1}
 =\frac{2\ell}{(2c-1)(2c-\ell-1)}.
\end{equation}
For $d=2k+1$ and $c=d/4+o(d)$, both denominators on the right are $k(1+o_k(1))$. Thus a fixed loss $\ell$ gives width $(2\ell+o_k(1))k^{-2}$. This interval uses the graph's actual $c$; estimating $c$ introduces a separate error.

\begin{lemma}[Uniform transfer]\label{s3:lem:transfer}\label{lem:fixedtargets}
For each row below and every fixed degree in that row, a.a.s. every balanced signing of $G_{n,d}$ with $c_\sigma\ge c_0$ satisfies
\begin{equation}\label{s3:eq:loss2}
 |\cut(S)|>q_0|\sigma(S)|
 \qquad(\varnothing\ne S\ne V(G)).
\end{equation}
In particular, $q_G(\sigma)>q_0$.
\[
\begin{array}{c|c|c}
 d&c_0&q_0\\\hline
 d\ge21&d/4+(19/60)\sqrt d&2c_\sigma-2\\
 15&499/100&801/100\\
 18&6&10\\
 19&613/100&101/10
\end{array}
\]
\end{lemma}

\begin{proof}[Proof of Theorems~\ref{main:jaeger} and~\ref{main:odd-flow}, and of~\eqref{eq:loss2intro}]
Let $d\ge21$ be odd, and choose a maximum-bisection signing $\sigma$, of density $c$. Lemma~\ref{s1:lem:bisectionstrong} gives $c>d/4+\sqrt d/\pi$, so the uniform row of Lemma~\ref{s3:lem:transfer} applies to $\sigma$. Taking $\ell=2$ in~\eqref{s2:eq:window} gives \eqref{eq:loss2intro}. Since
\[
 2c-3>\frac{d-6+(4/\pi)\sqrt d}{2}
 \qquad\text{and}\qquad
 \frac{4\sqrt d}{\pi}\ge\frac{4\sqrt{21}}{\pi}>5,
\]
this gives the last display of Theorem~\ref{main:odd-flow}. If $d=4p+1$, the same bounds give $2c-2>(d+1)/2=2p+1$. The strict cut inequalities imply that the graph is connected, so Lemma~\ref{s1:lem:criterion}(1) with $q=2p+1$ gives a strongly connected Jaeger orientation.

For odd $5\le d\le19$, Theorem~\ref{thm:p1strict} gives degree five, Propositions~\ref{prop:degree7} and~\ref{prop:degree11} give degrees seven and eleven, and Proposition~\ref{prop:small-strong} gives degrees nine, thirteen and seventeen. At degrees fifteen and nineteen, apply the corresponding rows of Lemma~\ref{s3:lem:transfer} to a maximum bisection. Lyons's estimate supplies their density thresholds. For degree nineteen the bound $19/4+\sqrt{19}/\pi>613/100$ suffices. For degree fifteen, use the exact bound and
\[
 \arcsin x>x+x^3/6+3x^5/40\quad(0<x<1),
 \qquad \sqrt{14}>187/50,\quad \pi<22/7,
\]
to obtain $c(G_{n,15})>499/100$ a.a.s. The resulting cut ratios exceed $801/100$ and $101/10$, respectively, and Lemma~\ref{s1:lem:criterion}(2) gives the displayed bounds. Finally, Theorem~\ref{thm:p1strict} and Proposition~\ref{prop:small-strong} supply strongly connected Jaeger orientations for $p\le4$.
\end{proof}

\section{The joint count and the uniform transfer}\label{sec:counting}

We prove Lemma~\ref{s3:lem:transfer} by counting a balanced signing together with a set that violates its cut inequality. We call such a set \emph{deficient}.

\subsection{The counting formula}

We continue in the pairing model introduced above. Conditional on simplicity the graph has distribution $G_{n,d}$, and for fixed $d$ the simplicity probability tends to \cite[Section~2.2, equation~(3)]{W}
\begin{equation}\label{s1:eq:simplicity}
 \exp(-(d^2-1)/4)>0.
\end{equation}
Thus a.a.s. statements in the pairing model pass to $G_{n,d}$.

Put $L(z)=z\log z$, with $L(0)=0$. We use the following elementary counting formula. Suppose the vertices are divided into a fixed number of ordered classes of proportions $\alpha_i$. Let $m_{ij}$ count edges between classes, with $m_{ii}$ counting internal edges, and put
\[
 P_{ii}=\frac{2m_{ii}}{dn},\qquad P_{ij}=\frac{m_{ij}}{dn}\quad(i\ne j).
\]
Then $P$ is symmetric, has row sums $\alpha_i$, and has total mass one. We call it the ordered-endpoint table. Each edge contributes its two ordered endpoint pairs, and $P_{ij}$ is the proportion with first vertex in class $i$ and second vertex in class $j$. The expected number of partitions with these class sizes and edge counts is
\[
 \frac{n!}{\prod_i(n\alpha_i)!}
 \frac{\prod_i(dn\alpha_i)!}{\prod_i2^{m_{ii}}m_{ii}!\prod_{i<j}m_{ij}!}
 \frac1{(dn-1)!!}.
\]
Here $(dn-1)!!=(dn-1)(dn-3)\cdots1$ is the total number of pairings. Define the counting exponent by
\begin{equation}\label{s1:eq:rate}
 \mathcal F_d(\alpha,P)=(d-1)\sum_iL(\alpha_i)-\frac d2\sum_{i,j}L(P_{ij}).
\end{equation}
The estimate $\log(m!)=m\log m-m+O(\log(m+2))$, obtained by comparing $\sum_{j=1}^m\log j$ with its integral, shows that the logarithm of the displayed expected count is $n\mathcal F_d(\alpha,P)+O_d(\log n)$. The error is uniform, including zero entries. Since there are polynomially many patterns, a uniform negative bound on $\mathcal F_d$ excludes all such partitions a.a.s.

Counting signings as well as sets introduces an exponential factor that is too costly for very small sets. We first bound small cuts independently of the signing. The following criterion will also be used in Appendix~\ref{s3:app:small}.

For $0<\theta<1$ and $0\le a<1/(1+\theta)$, define
\begin{equation}\label{s1:eq:smallrate}
\begin{split}
 F_{d,\theta}(a)={}&(d-1)\{L(a)+L(1-a)\}\\
 &-\frac d2\{L((1-\theta)a)+2L(\theta a)
                 +L(1-(1+\theta)a)\}.
\end{split}
\end{equation}
To compare degrees, write $H(a)=-L(a)-L(1-a)$ for $a>0$. For the two-class table
$P=\bigl(\begin{smallmatrix}(1-\theta)a&\theta a\\ \theta a&1-(1+\theta)a\end{smallmatrix}\bigr)$ with row sums $\alpha=(a,1-a)$, define its relative entropy against the product table by
\[
 I(P)=\sum_{i,j}P_{ij}\log\frac{P_{ij}}{\alpha_i\alpha_j},
\]
with zero summands omitted. The inequality $u\log(u/v)\ge u-v$ shows that $I(P)\ge0$, and $F_{d,\theta}(a)=H(a)-dI(P)/2$ therefore decreases with $d$. Differentiation also gives
\begin{equation}\label{s1:eq:smallsecond}
 F_{d,\theta}''(a)=
 \frac{d(1-\theta)-2-(d-2)(1+\theta)a}
 {2a(1-a)(1-(1+\theta)a)}.
\end{equation}

\begin{lemma}\label{s1:lem:smallcriterion}
Fix $d$, $0<\theta<1$ and $0<\eta<1-\theta$, with $d(1-\theta)>2$. Suppose that, for some $\epsilon>0$,
\[
 F_{d,\theta}(a)\le-\epsilon a\qquad(0<a\le\eta).
\]
Then a.a.s. every $S$ with $0<|S|\le\eta n$ has
$|\cut(S)|>\theta d|S|$.
\end{lemma}

\begin{proof}
Put $a=|S|/n$ and $b_{\rm cut}=|\cut(S)|/(dn)$. The two-class table is
\[
 P=\begin{pmatrix}a-b_{\rm cut}&b_{\rm cut}\\ b_{\rm cut}&1-a-b_{\rm cut}\end{pmatrix}.
\]
The derivative of \eqref{s1:eq:rate} with respect to $b_{\rm cut}$ is
\[
 \frac d2\log\frac{(a-b_{\rm cut})(1-a-b_{\rm cut})}{b_{\rm cut}^2}>0
 \qquad(0<b_{\rm cut}<a(1-a)).
\]
Since $\theta<1-a$, all bad cut sizes $b_{\rm cut}\le\theta a$ have exponent at most $F_{d,\theta}(a)$. The uniform counting error in \eqref{s1:eq:rate} therefore excludes $C_{\rm size}\log n\le|S|\le\eta n$, for a sufficiently large fixed $C_{\rm size}$.

Put $\xi=d(1-\theta)/2>1$. A bad set of $s$ vertices contains at least $m_s=\lceil\xi s\rceil$ internal pairs. If $2m_s>ds$ it is impossible. Otherwise choosing these points and pairing them gives
\[
 \binom ns\binom{ds}{2m_s}\frac{(2m_s-1)!!}{(dn-2m_s)^{m_s}}
 \le\left[C_d(s/n)^{\xi-1}\right]^s
 \qquad(1\le s<C_{\rm size}\log n),
\]
for a constant $C_d$ and all sufficiently large $n$. Indeed, $m_s\ge\xi s$ and $m_s=O_d(s)$, and the pairing probability is at most $(C_ds/n)^{m_s}$. Since $\xi>1$, summing gives $o(1)$. Conditioning on simplicity completes the proof.
\end{proof}

\subsection{Comparing the two cuts}

For a balanced signing with positive side $A$ and a set $S$, the four ordered classes are $A\cap S$, $A^c\cap S$, $A\setminus S$, $A^c\setminus S$, with proportions
\[
 \alpha=(a,b,x,y),\qquad x=1/2-a,\quad y=1/2-b.
\]
After complementing $S$ and reversing signs, their parameters belong to
\[
 \mathcal D_0=\{(a,b):a\ge b\ge0,\quad0\le a+b\le1/2\}.
\]
Let $\chi_{ij}$ be one when classes $i,j$ have different signs and zero otherwise. Let $\psi_{ij}$ be one when exactly one of the two classes is in $S$ and zero otherwise. The class arrangement is
\[
 \begin{array}{c|cc}
       &S&S^c\\\hline
 A     &a&x\\
 A^c   &b&y
 \end{array}
 \qquad\text{with endpoint-table class order }(a,b,x,y).
\]
The same table and its two cut statistics will be used in the later transfer arguments.

\begin{lemma}\label{s1:lem:comparison}
Let $P$ be the four-class endpoint table with row sums $\alpha=(a,b,x,y)$, omitting any empty classes and their zero rows. Put
\[
 c_\sigma=\frac d2\sum_{i,j}\chi_{ij}P_{ij},\qquad
 \zeta=\frac d2\sum_{i,j}\psi_{ij}P_{ij}.
\]
These are the bisection density and $|\cut(S)|/n$, respectively. Choose positive class weights $v_i$ and positive cut multipliers $U,V$, and define
\[
 R_{ij}=v_iv_jU^{\chi_{ij}}V^{\psi_{ij}},\qquad W=\sum_{i,j}R_{ij}.
\]
Then
\begin{equation}\label{s1:eq:comparison}
 \mathcal F_d\le(d-1)\sum_iL(\alpha_i)
       -d\sum_i\alpha_i\log v_i
       -c_\sigma\log U-\zeta\log V+\frac d2(W-1).
\end{equation}
\end{lemma}

\begin{proof}
The inequality $u\log(u/v)\ge u-v$ gives
\[
 -\sum P_{ij}\log P_{ij}
 \le-\sum P_{ij}\log R_{ij}+W-1.
\]
The row weights contribute $-2\sum_i\alpha_i\log v_i$ to the last sum. The cut multipliers contribute $-(2c_\sigma/d)\log U-(2\zeta/d)\log V$. Substitution into \eqref{s1:eq:rate} proves the bound.
\end{proof}

\subsection{Proof of the transfer}

We also use the following bound on all bisections; its two-class proof is in Appendix~\ref{s3:app:scalar}.

\begin{lemma}\label{s3:lem:scalar}
For every fixed $d\ge3$, a.a.s. every bisection has density less than $293d/1000+1$.
\end{lemma}

\begin{proof}[Proof of the uniform row of Lemma~\ref{s3:lem:transfer}]
Write $c=c_\sigma$ and $A=\sigma^{-1}(1)$, and suppose $S$ is deficient. Replacing $S$ by its complement does not change its cut size or $|\sigma(S)|$, and replacing $\sigma$ by $-\sigma$ preserves the bisection density. We may therefore assume $|S|\le n/2$ and $\sigma(S)\ge0$. Write $A^c=V(G)\setminus A$, and put
\[
 a=|A\cap S|/n,\quad b=|A^c\cap S|/n,\quad
 x=1/2-a,\quad y=1/2-b,
 \qquad t=x+b=|A\mathbin\triangle S|/n.
\]
Thus $t$ is the relative size of the symmetric difference of $S$ and the positive bisection side. Also $a\ge b\ge0$, $a+b\le1/2$, and $t=1/2-(a-b)$. We treat the following three regions.
\begin{center}\small
\begin{tabular}{@{}lll@{}}\toprule
Region&Method&Proof\\\midrule
$|S|\le7n/50$&two-class count, before $\sigma$ is chosen&Appendix~\ref{s3:app:small}\\
$t\le2/9$&quadratic entropy bound near the side $A$&Appendix~\ref{s3:app:near}\\
$(a,b)$ in the pentagon \eqref{s3:eq:domain}&four-class comparison and endpoint lemma&Appendices~\ref{s3:app:global}, \ref{s3:app:tables}\\\bottomrule
\end{tabular}
\end{center}

For $d\ge21$, small sets, with $|S|\le7n/50$, are excluded by the small-set estimate in Appendix~\ref{s3:app:small} and Lemma~\ref{s3:lem:scalar}:
\[
 |\cut(S)|>\frac{293d}{500}|S|
 >(2c-2)|\sigma(S)|
\]
when $\sigma(S)\ne0$; if $\sigma(S)=0$, the positive first bound already suffices. Next, the quadratic entropy bound in Appendix~\ref{s3:app:near} excludes every deficient signing--set pair with $t\le2/9$.

The remaining pairs have $(a,b)$ in the fixed pentagon
\begin{equation}\label{s3:eq:domain}
 \mathcal D_\eta=\{a\ge b\ge0:\ \eta\le a+b\le1/2,\quad a-b\le5/18\},
 \qquad \eta=7/50.
\end{equation}
Lemma~\ref{lem:endpoint-reduction} in Appendix~\ref{s3:app:global} bounds the counting exponent $\mathcal F_d$ by a function of $t$. The endpoint estimates in Appendix~\ref{s3:app:tables} make this bound uniformly negative throughout the pentagon.

In both regions treated by entropy, each class-size and edge-count pattern has expected count at most $\exp(-\epsilon n+O_d(\log n))$ for some fixed $\epsilon>0$. Summing over the polynomially many patterns gives $o(1)$. Together with the small-set estimate, this excludes all deficient pairs in the pairing model. Conditioning on simplicity, using~\eqref{s1:eq:simplicity}, completes the proof.
\end{proof}

\begin{proof}[Proof of the finite-degree rows of Lemma~\ref{s3:lem:transfer}]
Let $d\in\{15,18,19\}$, and let $c_0,q$ be the corresponding density and cut thresholds in the lemma. Put $\eta=13/100$. Appendix~\ref{s3:app:small} excludes all sets with $0<|S|\le\eta n$.

For any remaining deficient set, complementing and reversing signs gives the four-class coordinates above and $h=a-b=1/2-t$. Changing membership of $tn$ vertices changes a cut by at most $dtn$. Thus
\[
 c_0\le |\cut(S)|/n+dt\le q/2+(d-q)t,
 \qquad t\ge t_0:=\frac{c_0-q/2}{d-q}.
\]
It remains to count pairs on the pentagon
\[
 \mathcal P_{d}=\{a\ge b\ge0:\ \eta\le a+b\le1/2,\quad a-b\le1/2-t_0\}.
\]
The one-variable bounds in Appendix~\ref{s3:app:tables} give a comparison bound below $-1/200$ throughout each pentagon. Summing over the polynomially many patterns excludes all remaining pairs in the pairing model. Conditioning on simplicity finishes the proof.
\end{proof}

\section{The one-unit gap between bisections and orientations}\label{s3:sec:gap}

Since $q=d-2p$, a loss of two in the cut parameter costs one in the orientation parameter. We apply this observation to the uniform transfer for $d\ge21$ and then treat the smaller degrees. Throughout, $G=G_{n,d}$.

\subsection{The bisection obstruction}

The observation in Section~\ref{sec:introduction} gives $p_{\rm bis}(d)\le p^*(d)$ whenever $p^*(d)<d/2$. To include the case $p^*(d)=d/2$, choose a bisection uniformly at random. Each edge is separated with probability $n/(2(n-1))$, so some bisection has at least $dn/4$ edges. Consequently the deterministic averaging bound gives
\begin{equation}\label{s3:eq:meanbound}
 p_{\rm bis}(d)\le\lceil d/4\rceil<d/2\qquad(d\ge3).
\end{equation}

\subsection{Augmentation and the upper bound}

\begin{lemma}[Adding a perfect matching]\label{lem:matching-augmentation}
Fix an integer $d\ge3$ and a real number $r\ge0$.
If $q^*(G_{n,d})>r$ a.a.s., then $q^*(G_{n,d+1})>r$ a.a.s.
\end{lemma}
\begin{proof}
Let $G$ have distribution $G_{n,d}$, and let $M$ be an independent
uniform perfect matching of its vertices. For every balanced signing
$\sigma$ and every vertex set $S$,
\[
 |\cut_{G\cup M}(S)|\ge |\cut_G(S)|,
 \qquad q_{G\cup M}(\sigma)\ge q_G(\sigma).
\]
The probability that the union is simple is bounded away from zero
\cite[proof of Lemma~4.14]{W}, so conditioning on this event preserves
the asserted a.a.s. property. The conditioned union and $G_{n,d+1}$
are mutually contiguous by \cite[Theorem~4.15(ii)]{W}.
The conclusion follows.
\end{proof}

\begin{proof}[Proof of Theorem~\ref{main:gap}]
The deterministic lower bound was proved above. For $d\ge21$, we first prove \eqref{eq:pointwiseintro}.
Put $c=\bis(G)/n$ and $p_0=\lceil d/2-c\rceil$. The bisection obstruction gives the lower bound for any orientation parameter less than $d/2$; if the parameter equals $d/2$, the same bound follows from $p_0\le\lceil d/4\rceil<d/2$. Lemma~\ref{s3:lem:scalar} ensures $c<d/2$ and hence $p_0\ge1$. Averaging gives $c\ge d/4$, so
\[
 p_0+1\le\lceil d/4\rceil+1<d/2.
\]
Choose a maximum bisection. It satisfies the density hypothesis by Lemma~\ref{s1:lem:bisectionstrong}. Apply the uniform row of Lemma~\ref{s3:lem:transfer} and Lemma~\ref{s1:lem:criterion}(1) with the positive integer $q=d-2(p_0+1)\le2c-2$, which has the parity of $d$. This constructs the orientation required for the upper bound.

If $p=p_{\rm bis}(d)$, then a.a.s. $c\ge d/2-p$ and hence $p_0\le p$. The same cut inequality, now with the fixed parameter $q=d-2(p+1)\le d-2(p_0+1)$, gives a $(p+1)$-orientation by Lemma~\ref{s1:lem:criterion}(1). Hence $p^*(d)\le p_{\rm bis}(d)+1$.

For the smaller degrees, adding an independent Hamilton cycle and
directing it cyclically gives
\begin{equation}\label{eq:augmentation}
 G_{n,d}\text{ a.a.s. has a $p$-orientation}
 \quad\Longrightarrow\quad
 G_{n,d+2}\text{ a.a.s. has a $(p+1)$-orientation}.
\end{equation}
Indeed, every outdegree increases by one. Conditioning on a simple
union has probability bounded away from zero, and the conditioned
union is mutually contiguous to $G_{n,d+2}$
\cite[Theorem~4.15(i) and proof of Lemma~4.14]{W}.

For each even $10\le d\le20$, the cut-ratio bounds used in the proof of Theorem~\ref{main:odd-flow} give $q^*(G_{n,d-1})>d/2$ a.a.s. Lemma~\ref{lem:matching-augmentation} preserves this inequality at degree $d$. Taking the even divergence parameter $q=2\lfloor d/4\rfloor$ in Lemma~\ref{s1:lem:criterion}(1) yields
\begin{equation}\label{eq:small-even-orientations}
 p^*(d)\le\lceil d/4\rceil
 \qquad(d\in\{10,12,14,16,18,20\}).
\end{equation}

The following table records all the inputs. The scalar bound in
Lemma~\ref{s3:lem:scalar} gives $p_{\rm bis}\ge2$ at $d\ge10$ and
$p_{\rm bis}\ge3$ at $d\ge15$. It also gives $p_{\rm bis}(20)\ge4$. At $d=14,19$ the stronger lower
bounds follow from the maximum-cut estimates in Appendix~\ref{s3:app:clms}.
\begin{center}\small
\setlength{\tabcolsep}{5pt}
\begin{tabular}{@{}clcc@{}}\toprule
$d$&Orientation input&Lower bound on $p_{\rm bis}$&Upper bound on $p^*$\\\midrule
$3,4$&Hamiltonicity \cite[Theorem~1]{RW}&$1$&$1$\\
$5,6$&Theorem~\ref{thm:p1strict}; augment $4$&$1$&$1,2$\\
$7,8,9$&Augment $5$; add a matching to $7$; \cite[Theorem~1.3]{LLL}&$1$&$2$\\
$10,11,12$&\eqref{eq:small-even-orientations}; augment $9$&$2$&$3$\\
$13,14$&\cite[Theorem~1.3]{LLL}; \eqref{eq:small-even-orientations}&$2,3$&$3,4$\\
$15,16,17$&Augment $13$; \eqref{eq:small-even-orientations}; \cite[Theorem~1.3]{LLL}&$3$&$4$\\
$18,19,20$&\eqref{eq:small-even-orientations}; augment $17$&$3,4,4$&$5,5,5$\\\bottomrule
\end{tabular}\end{center}
At degree three, direct a Hamilton cycle cyclically and the remaining
matching arbitrarily. At degree four, direct the even Hamilton cycle
alternately and each cycle of the remaining $2$-regular graph
cyclically; the outdegrees are one or three. Thus the first row uses
the Hamiltonicity theorem. At degree eight, Proposition~\ref{prop:degree7} and Lemma~\ref{lem:matching-augmentation} give $q^*(G_{n,8})>21/5$ a.a.s. Lemma~\ref{s1:lem:criterion}(1) with $q=4$ then gives a $2$-orientation.

Each row gives $p^*\le p_{\rm bis}+1$, except possibly $d=18$
when $p_{\rm bis}=3$. In that case its definition supplies a bisection
of density at least six. The degree-eighteen row of
Lemma~\ref{s3:lem:transfer} and Lemma~\ref{s1:lem:criterion}(1),
with $q=10$, give a $4$-orientation. This completes the proof.
\end{proof}

\section{Smaller losses at large degree}\label{s2:sec:precision}

Retaining the bisection density in the joint count gives losses below two. We use these estimates to sharpen the flow interval and the rounding bound for orientations.

The sharper asymptotic loss uses the Gaussian maximum-bisection constant. To specify its normalization, for each integer $N\ge2$ let $g_{ij}$, $1\le i<j\le N$, be independent Gaussian random variables of mean zero and variance one. Put
\begin{equation}\label{eq:SKconstant}
 C^*=\lim_{N\to\infty}\frac1N\mathbb E
       \max_{\tau\in\{-1,1\}^N}
       \frac1{\sqrt N}\sum_{1\le i<j\le N}g_{ij}\tau_i\tau_j.
\end{equation}
The constant is denoted by $P_*$ in \cite[Section~1.1, equations~(1.6)--(1.7)]{DMS}, and its numerical estimate is $C^*\approx0.76316673$ \cite[after equation~(3)]{OSS}. The bisection correction is $(C^*/2)\sqrt d$ \cite[Theorem~1.5]{DMS}.

\subsection{A transfer estimate which retains the bisection density}

Li and Li \cite[Theorem~3.7]{LL} give a signing with cut parameter $d/2+(C^*-\epsilon)\sqrt d$ for each fixed $\epsilon>0$ at sufficiently large odd degree. The next theorem measures the loss from $2c_\sigma$, uniformly over all signings of the stated density. The density constraint reduces the counting exponent and allows a smaller loss.

Let $z_0>1$ be the unique root of $z^3-z^2-z-1=0$, and put $T_0=z_0^2$. The root $z_0$ is the tribonacci constant. Uniqueness follows because the polynomial is negative at one, tends to infinity, and has derivative $(3z+1)(z-1)>0$ for $z>1$.

\begin{theorem}\label{s2:thm:transfer}
Fix $\beta>0$ and
\begin{equation}\label{s2:eq:losscondition}
 \ell>\max\left\{0,\frac{2(\log2-4\beta^2)}{\log T_0}\right\}.
\end{equation}
There is $d_0=d_0(\beta,\ell)$ such that, for every fixed $d\ge d_0$, a.a.s. every balanced signing with
$c_\sigma\ge d/4+\beta\sqrt d$ satisfies
\[
 |\cut(S)|>(2c_\sigma-\ell)|\sigma(S)|
 \qquad(\varnothing\ne S\ne V(G_{n,d})).
\]
\end{theorem}

\begin{proof}
Write $c=c_\sigma$ within this proof. Call a nontrivial set deficient if its cut is at most $(2c-\ell)|\sigma(S)|$. Complementing $S$ and reversing all signs allows us to assume $|S|\le n/2$ and $\sigma(S)\ge0$. Put
\[
 u=\frac12-\frac{|S|}n,\qquad
 t=\frac12-\frac{\sigma(S)}n.
\]
Then $0\le u\le t\le1/2$. If $A=\sigma^{-1}(1)$, then $tn$ is the number of vertices belonging to exactly one of $S,A$. We first show that a deficient set must differ from $A$ at few vertices.

\paragraph{Localization.} Let $M$ be the adjacency matrix, and let $\lambda$ be the largest absolute value of its eigenvalues on the subspace perpendicular to the all-ones vector. Friedman's theorem, in the form \cite[Lemma~3.2]{LL}, gives $\lambda\le2\sqrt d<d$ a.a.s. for every fixed $d\ge5$. The spectral cut bound \cite[Lemma~3.4]{LL}, together with the quadratic form for a balanced bisection, gives
\[
 c\le\frac{d+\lambda}{4},\qquad
 \frac{|\cut(S)|}{n}\ge(d-\lambda)\left(\frac14-u^2\right).
\]
The case $t=1/2$ would give $\sigma(S)=0$ and $|\cut(S)|=0$, contradicting $\lambda<d$. For $t<1/2$, deficiency and $u\le t$ give
\[
 (d-\lambda)(1/2-t)(1/2+t)
 \le(2c-\ell)(1/2-t)
 \le\left(\frac{d+\lambda}{2}-\ell\right)(1/2-t).
\]
Cancel the positive factor $1/2-t$ to obtain
\begin{equation}\label{s2:eq:localization}
 (d-\lambda)t\le\lambda-\ell,
 \qquad t\le\frac{\lambda}{d-\lambda}=O(d^{-1/2}).
\end{equation}
The case $t=0$ is also impossible, since then $S=A$ and deficiency says $c\le c-\ell/2$. In particular, for any fixed $\eta>0$, all deficient sets lie in $0<t\le\eta$ once $d$ is sufficiently large.

\paragraph{The joint count.} We exclude this region by counting the signing and its deficient set together. Order the four classes as
\[
 A\cap S,\quad A^c\cap S,\quad A\setminus S,\quad A^c\setminus S,
\]
with proportions $\alpha=(a,b,x,y)$, where
\[
 a=1/2-x,\quad y=1/2-b,\quad b+x=t,\quad
 h=a-b=1/2-t.
\]
Use the ordered-endpoint table and the two indicator tables $\chi,\psi$ from Lemma~\ref{s1:lem:comparison}. The two cut constraints are
\begin{equation}\label{s2:eq:joint}
 \frac d2\sum_{i,j}\chi_{ij}P_{ij}=c,\qquad
 \frac d2\sum_{i,j}\psi_{ij}P_{ij}\le(2c-\ell)h.
\end{equation}
For fixed class sizes and edge counts, \eqref{s1:eq:rate} gives the logarithm of the expected count as $n\mathcal F_d+O_d(\log n)$. This count includes the choice of the signing.

Put $c_0=d/4+\beta\sqrt d$. The comparison separates cut weights from vertex weights. The factor $K$ is the odds ratio for the required crossing fraction $2c_0/d$, while $T>1$ penalizes a deficient cut. The symmetric weight $w_0$ normalizes the table when the two small classes are empty. A common weight $\kappa$ on these classes cancels their logarithmic terms; its choice below minimizes the limiting linear coefficient. Define
\[
 K=\frac{d+4\beta\sqrt d}{d-4\beta\sqrt d},\qquad
 w_0=[2(1+K)]^{-1/2}.
\]
For sufficiently large $d$ we have $d>16\beta^2$, so $K>1$. Choose $1<T<T_0$ later and put
\[
 \kappa=\frac2{T+T^{-1}},\qquad
 v=(w_0,\kappa b,\kappa x,w_0),\qquad
 R_{ij}=v_iv_j(KT)^{\chi_{ij}}T^{-\psi_{ij}},\qquad W=\sum_{i,j}R_{ij}.
\]
Apply Lemma~\ref{s1:lem:comparison} with multipliers $KT$ and $T^{-1}$, omitting zero rows as in that lemma. The row weights contribute $-d\sum_i\alpha_i\log v_i$. The two cut terms are at most
\[
 -c\log(KT)+(2c-\ell)h\log T
 =-c\log K-2ct\log T-\ell h\log T.
\]
This expression decreases with $c$, so replacing $c$ by $c_0$ gives an upper bound.

\paragraph{Two scalar bounds.} Direct summation gives
\begin{align*}
 W&=1+2w_0\kappa(KT+T^{-1})t
              +\kappa^2(b^2+x^2+2Kbx)\\
 &\le1+2w_0\kappa(KT+T^{-1})t
              +\frac{\kappa^2(K+1)}2t^2.
\end{align*}
Define $J=L(1/2-x)+L(1/2-b)$. For $0\le z\le\eta<1/2$ and any $\vartheta>0$, we have
\begin{align*}
 -L(z)&\le\vartheta-z(1+\log\vartheta),\\
 -\log2-(1-\log2)t
 &\le J\le-\log2-(1-\log2)t+\frac{t^2}{1-2\eta}.
\end{align*}
The first inequality is the tangent bound to the concave function $-L$ at $\vartheta$. For the bounds on $J$, use $L'(1/2)=1-\log2$ and $L''(z)=1/z\le2/(1-2\eta)$ for $1/2-\eta\le z\le1/2$, then use $b^2+x^2\le t^2$.

The small-class weights give the cancellation
\[
 (d-1)L(b)-db\log(\kappa b)=-L(b)-db\log\kappa,
\]
and likewise for $x$. Before applying the scalar bounds, the estimate is
\begin{align*}
 \mathcal F_d\le{}&(d-1)J-d(1-t)\log w_0-L(b)-L(x)-dt\log\kappa\\
 &-c_0\log K-2c_0t\log T-\ell(1/2-t)\log T+\frac d2(W-1).
\end{align*}
Use the upper bound on $dJ$ and the lower bound on $J$ in $(d-1)J=dJ-J$. Denote the constant term by $I_d$, the coefficient proportional to $dt$ by $A_d$, the remaining linear coefficient by $V_1$, and the quadratic bound by $B_d$. They are
\begin{align*}
 I_d&=\log2+\frac d2\log\frac{1+K}{2}-c_0\log K-\frac\ell2\log T,\\
 A_d&=-1+\log(2w_0)-\log\kappa-\tfrac12\log T
                         +w_0\kappa(KT+T^{-1}),\\
 V_1&=-\log2-\log\vartheta+\ell\log T,\qquad
 B_d=\frac1{1-2\eta}+\frac{\kappa^2(K+1)}4.
\end{align*}
The resulting bound is
\begin{equation}\label{s2:eq:entropybound}
 \mathcal F_d\le I_d+2\vartheta+
 t\{dA_d-2\beta\sqrt d\log T+V_1\}+dB_dt^2.
\end{equation}
The constant and linear terms explain the choice of parameters. For $0<\rho<1$, write
\[
 D_{\rm KL}(\rho\Vert1/2)=\rho\log(2\rho)+(1-\rho)\log(2(1-\rho))
\]
for its binary relative entropy against $1/2$. Set $\rho=K/(1+K)=1/2+2\beta/\sqrt d$. Then
\begin{align}
 I_d&=\log2-\frac d2D_{\rm KL}(\rho\Vert1/2)-\frac\ell2\log T\notag\\
 &\le\log2-4\beta^2-\frac\ell2\log T.
 \label{s2:eq:constant}
\end{align}
As a function of $\rho$, this relative entropy vanishes with its first derivative at $1/2$, and its second derivative is $1/[\rho(1-\rho)]\ge4$, giving $D_{\rm KL}(\rho\Vert1/2)\ge2(\rho-1/2)^2$. Moreover $K\to1$, $w_0\to1/2$, and
\[
 A_d\longrightarrow-\log(\kappa\sqrt T),
\]
because $\kappa(T+T^{-1})/2=1$. On writing $z=\sqrt T>1$, the inequality $\kappa\sqrt T>1$ becomes
\[
 z^4-2z^3+1=(z-1)(z^3-z^2-z-1)<0.
\]
It therefore holds exactly when $1<T<T_0$.

\paragraph{Choice of parameters.} By \eqref{s2:eq:losscondition}, choose such a $T$ with
\[
 m:=\frac\ell2\log T-\log2+4\beta^2>0,
 \qquad a_0:=\log(\kappa\sqrt T)>0.
\]
Choose $\vartheta>0$ so $2\vartheta<m/4$, and then $\eta>0$ small enough that
\[
 \eta\left\{\frac1{1-2\eta}+\frac{\kappa^2}{2}\right\}<a_0/4.
\]
For sufficiently large $d$, $A_d\le-3a_0/4$ and $\eta B_d\le a_0/4$. Dropping the negative square-root term in \eqref{s2:eq:entropybound} now yields, uniformly for $t\le\eta$,
\[
 \mathcal F_d\le-m+2\vartheta+t(-a_0d/2+V_1)<-m/2.
\]
For fixed $d$, the four class sizes and ten unordered edge counts give at most $(dn+1)^{14}$ patterns. Their total expected count is therefore at most a polynomial in $n$ times $e^{-mn/2}$, which tends to zero. Thus the pairing model a.a.s. has none of these pairs. The positive simplicity probability in \eqref{s1:eq:simplicity} transfers this assertion to $G_{n,d}$. On the simultaneous spectral event, \eqref{s2:eq:localization} excludes every other deficient set. This proves the theorem.
\end{proof}

\begin{corollary}\label{cor:explicit-loss-one}
For every fixed $d\ge2500$, a.a.s. every balanced signing with
$c_\sigma\ge d/4+(19/60)\sqrt d$ satisfies
\[
 |\cut(S)|>(2c_\sigma-1)|\sigma(S)|
 \qquad(\varnothing\ne S\ne V(G_{n,d})).
\]
Consequently, for every fixed odd $d\ge2501$, a.a.s., with $c=c(G_{n,d})$,
\[
 2+\frac2{2c-1}\le\phi(G_{n,d})<2+\frac2{2c-2}.
\]
\end{corollary}
\begin{proof}
Use the comparison in the proof of Theorem~\ref{s2:thm:transfer} with
\[
 \beta=\frac{19}{60},\quad \ell=1,\quad T=2,\quad
 \kappa=\frac45,\quad \eta=\frac1{24},\quad \vartheta=\frac1{1000}.
\]
Localization gives $t\le2/(\sqrt d-2)\le\eta$. The odds ratio and symmetric large-class weight satisfy
\[
 1<K\le\frac{769}{731},\qquad
 w_0=[2(1+K)]^{-1/2}\le\frac12.
\]
Thus the coefficients in \eqref{s2:eq:entropybound} obey
\[
 A_d\le\frac45(K-1)+\log5-\frac52\log2,\qquad
 B_d\le\frac{12}{11}+\frac4{25}(K+1).
\]
Consequently $A_d+\eta B_d<-1/50$. Moreover
\[
 I_d\le\tfrac12\log2-\frac{361}{900}<-\frac1{20},\qquad
 V_1=\log1000<7.
\]
Dropping the negative square-root term in \eqref{s2:eq:entropybound} therefore gives
\[
 \mathcal F_d<-\frac6{125}+t\left(-\frac d{50}+7\right)
 \le-\frac6{125}.
\]
Summing over the patterns and conditioning on simplicity complete the transfer. Lyons's bisection bound supplies the required density, and \eqref{s2:eq:window} gives the flow interval.
\end{proof}

For fixed $T$, the limiting linear coefficient for an arbitrary common small-class weight $\kappa>0$ is
\[
 -1-\log\kappa-\tfrac12\log T+\tfrac\kappa2(T+T^{-1}).
\]
Its derivative is $-1/\kappa+(T+T^{-1})/2$ and its second derivative is positive. Thus $\kappa=2/(T+T^{-1})$ is the unique minimizer among these weights, and the minimum is negative precisely for $1<T<T_0$.

For $\ell=2$, condition~\eqref{s2:eq:losscondition} holds for every $\beta>0$, since $T_0>2$. Thus a loss of two follows asymptotically from any positive gain of order $\sqrt d$ in the bisection density. Lemma~\ref{s3:lem:transfer} gives the explicit degree range.

Define
\[
 \ell_{\rm P}=\frac{2(\log2-C^{*2})}{\log T_0}.
\]
For any fixed $\ell>\max\{0,\ell_{\rm P}\}$, choose $\beta<C^*/2$ sufficiently close to $C^*/2$ that \eqref{s2:eq:losscondition} holds. The maximum-bisection asymptotic in \cite[Theorem~1.5]{DMS} gives $c(G_{n,d})\ge d/4+\beta\sqrt d$ a.a.s. for all sufficiently large fixed $d$. For large $d$ we also have $2c(G_{n,d})-\ell>1$. Theorem~\ref{s2:thm:transfer} and the cut criteria therefore give the following conclusion.

\begin{corollary}\label{s2:cor:window}
For every fixed $\ell>\max\{0,\ell_{\rm P}\}$ and every sufficiently large fixed odd $d$, a.a.s. \eqref{s2:eq:window} holds with $c=c(G_{n,d})$. Its width is given by \eqref{s2:eq:width}.
\end{corollary}

\begin{proof}[Proof of \eqref{eq:halfwindowintro}]
Take $\beta=19/60$ and $\ell=1/2$ in Theorem~\ref{s2:thm:transfer}. Condition~\eqref{s2:eq:losscondition} holds. Indeed, for $T=27/8$ we have $\kappa=2/(T+T^{-1})=432/793$ and $\kappa^2T>1$, so $T<T_0$ by the analysis of $\kappa\sqrt T$ in the proof of that theorem, and
\[
 \frac\ell2\log T-\log2+4\beta^2
 =\frac14\log\frac{27}8-\log2+\frac{361}{900}>\frac1{100}.
\]
Lemma~\ref{s1:lem:bisectionstrong} gives the required density for a maximum bisection. Substituting $\ell=1/2$ into \eqref{s2:eq:window} gives \eqref{eq:halfwindowintro}. Averaging gives $c\ge d/4$, while the spectral estimate gives $c\le(d+2\sqrt d)/4$. Hence both denominators in \eqref{s2:eq:width} are $k+O(\sqrt k)$ when $d=2k+1$, and the width is $(1+O(k^{-1/2}))k^{-2}$.
\end{proof}

\begin{corollary}\label{s2:cor:rounding}
Fix $\ell>\max\{0,\ell_{\rm P}\}$. For every sufficiently large fixed degree $d$, a.a.s., with $c=c(G_{n,d})$,
\[
 \left\lceil\frac d2-c\right\rceil
 \le p_{\rm or}(G_{n,d})
 \le\left\lceil\frac d2-c+\frac\ell2\right\rceil.
\]
\end{corollary}

\begin{proof}
Choose $\beta<C^*/2$ sufficiently close to $C^*/2$ for \eqref{s2:eq:losscondition}. The maximum-bisection theorem supplies density at least $d/4+\beta\sqrt d$, so Theorem~\ref{s2:thm:transfer} applies to a maximum-bisection signing. Put $p=\lceil d/2-c+\ell/2\rceil$. The spectral bound $c\le(d+2\sqrt d)/4$ and the averaging bound $c\ge d/4$ give $1\le p<d/2$ for large $d$. Thus $q=d-2p$ is a positive integer of the parity of $d$, and $q\le2c-\ell$. Lemma~\ref{s1:lem:criterion}(1) gives a $p$-orientation. Since this parameter is below $d/2$, the bisection obstruction applies to the minimum parameter and gives the lower bound.
\end{proof}

For $0<\ell<2$, the two ceilings agree whenever $j-1<d/2-c\le j-\ell/2$ for an integer $j$; in that case $p_{\rm or}(G)=j$.

Table~\ref{s2:tab:constants} lists the constants used in these bounds.
\begin{table}[ht]
\centering
\caption{Constants in the density-sensitive transfer. All numerical values are approximate.}\label{s2:tab:constants}
\begin{tabular}{@{}lll@{}}\toprule
Constant&Definition&Numerical value\\\midrule
$C^*$&Gaussian maximum in \eqref{eq:SKconstant}&$0.76316673$\\
$C^*/2$&Bisection correction coefficient&$0.38158336$\\
$2\sqrt2\,C^*$&Flow correction coefficient&$2.15856147$\\
$z_0$&$z_0^3-z_0^2-z_0-1=0$, $z_0>1$&$1.8392867552$\\
$T_0$&$z_0^2$&$3.3829757679$\\
$\ell_{\rm P}$&$2(\log2-C^{*2})/\log T_0$&$0.1816996$\\\bottomrule
\end{tabular}
\end{table}
The entries involving $C^*$ use the estimate reported in \cite[after equation~(3)]{OSS}.

\subsection{The limit of the joint first-moment count}

\begin{proposition}\label{prop:first-moment-barrier}
Fix $\beta,\ell>0$. For $d>4$, let $\mathcal R_d$ be the endpoint tables for balanced signings and sets with $\sigma(S)\ge0$, $|S|\le n/2$, satisfying
\[
 c_\sigma\ge c_d:=d/4+\beta\sqrt d,\qquad
 \frac{|\cut(S)|}{n}\le(2c_\sigma-\ell)(a-b),\qquad
 0<t=b+x\le\frac{2\sqrt d}{d-2\sqrt d}.
\]
The constraints are on real endpoint-table proportions with row sums $\alpha=(a,b,1/2-a,1/2-b)$, where $a\ge b\ge0$ and $a+b\le1/2$. Then
\[
 \lim_{d\to\infty}\sup_{(\alpha,P)\in\mathcal R_d}
 \mathcal F_d(\alpha,P)
 =\log2-4\beta^2-\frac\ell2\log T_0.
\]
In particular, when the right side is positive, this joint first-moment count cannot exclude all deficient pairs by a negative-exponent bound.
\end{proposition}
\begin{proof}
For the upper bound, fix $1<T<T_0$ and $\vartheta>0$. The estimates proving Theorem~\ref{s2:thm:transfer} give a negative coefficient of $dt$ once the fixed neighbourhood is sufficiently small and $d$ sufficiently large. Localization places all of $\mathcal R_d$ in that neighbourhood. Hence
\[
 \limsup_{d\to\infty}\sup_{\mathcal R_d}\mathcal F_d
 \le\log2-4\beta^2-\frac\ell2\log T+2\vartheta.
\]
Let $\vartheta\downarrow0$ and $T\uparrow T_0$.

For the matching lower bound, put
\[
 \rho_0=\frac{T_0^2}{1+T_0^2},\qquad
 s_\ell=\frac{\ell}{4\rho_0-3}
       =\ell\frac{1+T_0^2}{T_0^2-3},\qquad t_d=\frac{s_\ell}{d}.
\]
Since $T_0^2>3$, the number $s_\ell$ is positive, and $s_\ell>\ell$. Take
\[
 a=y=\frac{1-t_d}{2},\quad b=x=\frac{t_d}{2},\qquad
 D_d=\frac\ell2+(2c_d-\ell)t_d,
\]
and define
\[
 v_d=\frac{t_d+D_d/d}{4},\quad
 j_d=\frac{t_d-D_d/d}{4},\quad
 z_d=\frac{c_d}{d}-2v_d,\quad p_d=\frac12-t_d-z_d.
\]
For large $d$ these entries are nonnegative, and the symmetric table
\[
 P_d=\begin{pmatrix}
 p_d&v_d&j_d&z_d\\
 v_d&0&0&j_d\\
 j_d&0&0&v_d\\
 z_d&j_d&v_d&p_d
 \end{pmatrix}
\]
has row sums $(a,b,x,y)$. Its bisection density is $d(z_d+2v_d)=c_d$, while its set-cut density is
$d(z_d+2j_d)=c_d-D_d=(2c_d-\ell)(1/2-t_d)$. Also $t_d$ lies in the localization interval for all sufficiently large $d$.

The large entries have expansions
\[
 p_d=\frac14-\frac\beta{\sqrt d}-\frac{s_\ell-\ell}{4d}+O(d^{-3/2}),\qquad
 z_d=\frac14+\frac\beta{\sqrt d}-\frac{3s_\ell+\ell}{4d}+O(d^{-3/2}),
\]
and
\[
 dv_d\longrightarrow\frac{s_\ell\rho_0}{2},\qquad
 dj_d\longrightarrow\frac{s_\ell(1-\rho_0)}2.
\]
Substituting into \eqref{s1:eq:rate}, the $\log d$ terms from the two small classes cancel and Taylor expansion at $1/4$ gives
\[
 \mathcal F_d(\alpha,P_d)
 =\log2-4\beta^2-s_\ell D_{\rm KL}(\rho_0\Vert1/2)+o_d(1).
\]
The identity $T_0^2+1=2T_0^{3/2}$ follows from the defining tribonacci equation. Thus
\[
 \log(2\rho_0)=\tfrac12\log T_0,\qquad
 \log(2(1-\rho_0))=-\tfrac32\log T_0,
\]
so $s_\ell D_{\rm KL}(\rho_0\Vert1/2)=\ell\log T_0/2$. This proves the lower bound.
\end{proof}

An arbitrarily small increase of $v_d$ and decrease of $j_d$ preserves the row sums and makes both cut constraints strict, so rational and integer-count approximations give the same positive-exponent obstruction.

A positive exponent alone gives no lower bound on the probability of a deficient pair. The proposition identifies the limit of this joint first-moment bound.

\subsection{A quantitative form of the Li--Li asymptotic}

Li and Li \cite[Theorem~1.4]{LL} give an expansion with remainder $o_k(k^{-3/2})$. The interpolation estimate in \cite[Lemma~2.3, equation~(2.10)]{DMS} yields a remainder of $O(k^{-5/3})$. We first record the needed bisection estimate.

\begin{lemma}[Quantitative form of the bisection estimate in \cite{DMS}]\label{s2:lem:DMS}
There are absolute constants $B,d_0$ such that, for every fixed $d\ge d_0$,
\begin{equation}\label{s2:eq:DMS}
 \left|c(G_{n,d})-\frac d4-\frac{C^*}{2}\sqrt d\right|
 \le B d^{1/3}
\end{equation}
a.a.s. as $n\to\infty$.
\end{lemma}

\begin{proof}
Appendix~\ref{app:dms} derives the estimate from the dilute-graph comparison, whose error is $O(d^{1/3})$. The regularization and degree-shift errors are of order $d^{1/4}\sqrt{\log d}+(\log d)^2+\log d+1$. Conditioning on simplicity transfers the estimate from pairings to simple graphs.
\end{proof}

\begin{corollary}[Quantitative Li--Li asymptotic]\label{s2:thm:flowasymptotic}
There are absolute constants $K_0,K_{\rm flow}$ such that, for every fixed integer $k\ge K_0$, a.a.s.
\begin{equation}\label{s2:eq:flowformula}
 \phi(G_{n,2k+1})
 =2+\frac2k-\frac{2\sqrt2\,C^*}{k^{3/2}}+O(k^{-5/3}),
\end{equation}
where the absolute value of the remainder is at most $K_{\rm flow}k^{-5/3}$.
\end{corollary}

\begin{proof}
Put $d=2k+1$ and $H=2c(G_{n,d})-1$. Lemma~\ref{s2:lem:DMS} gives
\[
 H=k+\sqrt2\,C^*\sqrt k+O(k^{1/3}),\qquad
 \frac2H=\frac2k-\frac{2\sqrt2\,C^*}{k^{3/2}}+O(k^{-5/3}).
\]
The second estimate follows from $(1+x)^{-1}=1-x+x^2/(1+x)$ with $x=(H-k)/k=O(k^{-1/2})$. For any fixed admissible loss, \eqref{s2:eq:window} places the flow number between $2+2/H$ and $2+2/(H-\ell)$, whose difference is $O(k^{-2})$. This proves the claimed remainder.
\end{proof}

The $O(k^{-5/3})$ remainder comes from estimating $c(G_{n,d})$. The interval in Theorem~\ref{main:circular} uses $c(G_{n,d})$ itself and has length $O(k^{-2})$.

\section{A fixed gap at degree five}\label{sec:degreefive}

Pra\l{}at and Wormald \cite[Theorem~1.2]{PW} proved $\phi(G_{n,5})\le3$ a.a.s. We obtain a fixed gap below three by modifying a marked Jaeger orientation to give bounded monochromatic trees whose component graph satisfies weighted Hall inequalities.

\begin{theorem}\label{thm:p1strict}
With $\varepsilon_5=1/24006$, a.a.s. $G_{n,5}$ has a strongly connected Jaeger orientation and
\[
 \frac{14}{5}\le\phi(G_{n,5})<3-\varepsilon_5.
\]
\end{theorem}

\subsection{The forest obstruction}

Let $G$ be $5$-regular and let $\sigma$ be balanced with $q_G(\sigma)>3$. Each signing class induces a forest. Indeed, for every nonempty set $S$ within one class,
\[
 5|S|-2e(S)=|\cut(S)|>3|S|,
\]
so $e(S)<|S|$. A tree component of order $s$ has $3s+2$ boundary edges. Thus, if $C$ is the largest component order,
\begin{equation}\label{eq:p1-treebound}
 q_G(\sigma)\le3+\frac2C.
\end{equation}

\begin{proposition}\label{thm:p1forest}
Asymptotically almost surely, every partition of $G_{n,5}$ into two induced forests has more than $n/40$ vertices in tree components of order at least four. The partition need not be balanced.
\end{proposition}

Appendix~\ref{app:p1forest} proves the proposition by a weighted count of tree components. A coefficient bound and a one-variable maximization control the count for all class sizes.

\subsection{A matching criterion}

Fix a balanced partition $A\cup B$ into induced forests. Let $\mathcal A$ and $\mathcal B$ be their sets of tree components. Contract the trees to a bipartite multigraph $H$, retaining every crossing edge. For a tree $T$, write $s_T=|T|$, and for a collection $\mathcal X$ put $w(\mathcal X)=\sum_{T\in\mathcal X}s_T$. Since
\[
 \deg_H(T)=3s_T+2,
\]
the equality of the two crossing-degree sums and $|A|=|B|$ imply $|\mathcal A|=|\mathcal B|$. Put $C=\max_Ts_T$ and let $\sigma$ be positive on $A$.

A fractional perfect matching of $H$ assigns nonnegative weights to its edges with total weight one at each vertex. When $H$ has a perfect matching, define
\[
 a(H)=\max_z\min_{e\in E(H)}z_e,
\]
where $z$ ranges over fractional perfect matchings. The maximum exists by compactness, and $a(H)\le1/(3C+2)$.

\begin{proposition}\label{prop:p1matching}
For $0\le\delta\le2/C$, the inequality $q_G(\sigma)\ge3+\delta$ is equivalent to the existence of nonnegative edge weights $x_e$ on $H$ satisfying
\begin{equation}\label{eq:p1demands}
 \sum_{e\ni T}x_e=1-\frac\delta2s_T
 \qquad(T\in\mathcal A\cup\mathcal B).
\end{equation}
Equivalently, writing $N_H(\mathcal X)$ for the neighbours of $\mathcal X\subseteq\mathcal B$,
\begin{equation}\label{eq:p1weightedhall}
 |N_H(\mathcal X)|-|\mathcal X|
 \ge\frac\delta2\{w(N_H(\mathcal X))-w(\mathcal X)\}
 \qquad(\mathcal X\subseteq\mathcal B).
\end{equation}
In particular, $q_G(\sigma)\ge3$ if and only if $H$ has a perfect matching. In that case,
\begin{equation}\label{eq:p1matchingformula}
 q_G(\sigma)=\frac3{1-2a(H)},\qquad
 \phi(G)\le\frac3{1+a(H)}.
\end{equation}
\end{proposition}

\begin{proof}
A fractional orientation assigns to the two directions of each edge nonnegative weights summing to one. The real version of the prescribed-divergence criterion, obtained from the circulation criterion \cite[Theorem~1]{DGKNS}, realizes a zero-sum vector $b$ exactly when $|b(S)|\le|\cut(S)|$ for every $S$.

Apply this with $b=(3+\delta)\sigma$. Vertices of $A$ have fractional indegree $1-\delta/2$, and vertices of $B$ have fractional outdegree $1-\delta/2$. Let $x_e$ be the weight from $B$ to $A$ on a crossing edge. A tree of order $s_T$ has $s_T-1$ internal edges, so summing the prescribed degrees over it gives \eqref{eq:p1demands}.

Conversely, suppose \eqref{eq:p1demands} holds. Every $x_e$ is at most one. Assign that weight from $B$ to $A$ on a crossing edge and the remaining weight in the opposite direction. At $v$ in an $A$-tree $T$, let $x_v$ be the incoming crossing weight. Its required internal indegree is $z_v=1-\delta/2-x_v$. Since $x_v\le1-\delta s_T/2$, each $z_v$ is in $[0,1]$. Their sum over $T$ is $s_T-1$. Orient the tree away from a root chosen with probability $1-z_v$ at $v$. The average orientation has the prescribed internal indegrees. Reverse the construction on the $B$-trees. This realizes the prescribed divergence. The two sides of $H$ have equal total demand; applying the bipartite flow criterion gives \eqref{eq:p1weightedhall}. At $\delta=0$ this is Hall's criterion for a perfect matching.

Finally, the transformations
\[
 z_e=\frac{x_e+\delta/6}{1+\delta/3},
 \qquad x_e=(1+\delta/3)z_e-\delta/6
\]
convert \eqref{eq:p1demands} into a fractional perfect matching with $z_e\ge\delta/(6+2\delta)$. The identity $\deg_H(T)=3s_T+2$ verifies the vertex sums. Hence $q_G(\sigma)\ge3+\delta$ is equivalent to $a(H)\ge\delta/(6+2\delta)$. The bound $a(H)\le1/(3C+2)$ ensures that the full feasible range corresponds to $0\le\delta\le2/C$. Maximizing $\delta$ and applying Lemma~\ref{s1:lem:criterion}(2) gives \eqref{eq:p1matchingformula}.
\end{proof}

Two sufficient conditions follow. If $H$ has a perfect matching pairing trees of equal order, put $x_e=1-s_T/C$ on each matched edge and zero on the others. Then
\begin{equation}\label{eq:p1equalorders}
 q_G(\sigma)=3+2/C,
 \qquad \phi(G)\le3-\frac1{C+1}.
\end{equation}
This applies when all monochromatic trees have the same order $C$, since the component graph is $(3C+2)$-regular bipartite and has a perfect matching. Alternatively, if $K$ perfect matchings together cover every edge of $H$, their average gives $a(H)\ge1/K$, and
\[
 \phi(G)\le\frac{3K}{K+1}.
\]

For $0<\varepsilon<1$, a $5$-regular graph satisfies $\phi(G)\le3-\varepsilon$ if and only if it has a balanced two-forest partition whose component graph has a perfect matching and
\begin{equation}\label{eq:p1fixedtarget}
 a(H)\ge\frac{\varepsilon}{3-\varepsilon}.
\end{equation}
The reverse direction follows from \eqref{eq:p1matchingformula}; the forward direction follows from Lemma~\ref{s1:lem:criterion}(3) and the forest observation. Appendix~\ref{app:p1strict} obtains a lower bound on $a(H)$ independent of $n$ for the repaired partition.

The matching reduction holds for every odd $d\ge3$. For a balanced two-forest partition of a $d$-regular graph, replace $3s_T+2$ by $(d-2)s_T+2$ in the proof. If its component graph has a perfect matching, then
\[
 q_G(\sigma)=\frac{d-2}{1-2a(H)},\qquad
 \phi(G)\le\frac{2(d-2)}{d-3+2a(H)},
\]
with the flow bound understood when the denominator is positive.

\begin{proof}[Proof of Theorem~\ref{thm:p1strict}]
Put $\delta_5=1/12002$. Proposition~\ref{prop:p1repair} in Appendix~\ref{app:p1strict} supplies connectedness and a balanced signing with $q_G(\sigma)\ge3+\delta_5$ a.a.s. Lemma~\ref{s1:lem:criterion} gives the strong orientation and
\[
 \phi(G)\le 2+\frac{2}{2+\delta_5}
 =3-\frac{\delta_5}{2+\delta_5}<3-1/24006.
\]
For the lower bound, an optimal signing induces two forests because $q^*(G)>3$. Proposition~\ref{thm:p1forest} gives a tree of order at least four, so \eqref{eq:p1-treebound} gives $q^*(G)\le7/2$. The flow formula in Lemma~\ref{s1:lem:criterion}(3) yields $\phi(G)\ge14/5$.
\end{proof}

\section{Convergence at fixed degree}\label{sec:questions}

We ask whether the circular flow number converges in probability to a constant at each fixed degree.

\begin{conjecture}\label{conj:flow-limit}
For every fixed integer $d\ge3$ there exists a finite constant $f_d$ such that, for every $\varepsilon>0$,
\[
 \Pr\bigl(\lvert\phi(G_{n,d})-f_d\rvert>\varepsilon\bigr)
 \longrightarrow 0
 \qquad(n\longrightarrow\infty),
\]
where $n>d$ and $nd$ is even.
\end{conjecture}

For even $d$ we have $f_d=2$, since a unit flow in an Eulerian orientation gives a circular $2$-flow. The conjecture also holds at $d=3$, with $f_3=4$. In fact, $G_{n,3}$ is Hamiltonian a.a.s. by Robinson and Wormald \cite[Theorem~1]{RW}. An even Hamilton cycle splits into two perfect matchings $M_1,M_2$, and the remaining edges form a third perfect matching $M_3$. Orient the cycles in $M_1\cup M_2$ and in $M_2\cup M_3$ to give unit circulations. In a common reference orientation, add the first circulation to twice the second. On $M_2$ the resulting absolute values are $1$ or $3$, and on the other matchings they are $1$ or $2$. Reversing edges with negative values gives a circular $4$-flow.

A simple cubic graph that is not bipartite has a same-sign edge in every balanced signing. Its two endpoints give a cut of size $4$ and imbalance $2$, so $q^*(G)\le2$ and Lemma~\ref{s1:lem:criterion}(3) gives $\phi(G)\ge4$. To see that $G_{n,3}$ is not bipartite a.a.s., work first in the pairing model. Every bipartition has two equal classes, by regularity. For a fixed such partition, putting $m=3n/2$, the probability that all pairs cross is
\[
 \frac{m!}{(2m-1)!!}=\frac{2^m}{\binom{2m}{m}}\le(2m+1)2^{-m}.
\]
There are at most $2^n$ partitions to try, giving an upper bound of $(3n+1)2^{-n/2}$. Conditioning on simplicity preserves the $o(1)$ bound because its probability is bounded away from zero. Thus $\phi(G_{n,3})=4$ a.a.s.

\appendix
\section{Estimates for the bisection--orientation transfer}\label{s3:app:gap}

The first three subsections supply the bisection and small-cut estimates used in Sections~\ref{sec:counting} and~\ref{s3:sec:gap}. The remaining subsections develop the four-class comparison and its endpoint bounds.

\subsection{An upper bound for bisections}\label{s3:app:scalar}

The two-class count gives an upper bound on every bisection, which in turn fixes the small-cut threshold in the uniform transfer.

\begin{proof}[Proof of Lemma~\ref{s3:lem:scalar}]
Apply \eqref{s1:eq:rate} to the two classes of a bisection, with $\alpha=(1/2,1/2)$ and $P=\tfrac12\bigl(\begin{smallmatrix}1-r&r\\r&1-r\end{smallmatrix}\bigr)$. Here $r=2c/d$ is the fraction of all edges crossing the bisection. The counting exponent is
\[
 \log2-\frac d2D_{\rm KL}(r\Vert1/2),\qquad
 D_{\rm KL}(r\Vert1/2)=r\log(2r)+(1-r)\log(2(1-r)).
\]
Nonnegativity of the relative entropy of $(r,1-r)$ with respect to $(2/3,1/3)$ gives
\[
 D_{\rm KL}(r\Vert1/2)\ge r\log2-\log(3/2).
\]
At $r\ge293/500+2/d$, the counting exponent is at most
\[
 -\frac d2\{(793/500)\log2-\log3\}<0.
\]
Summing over cut sizes and conditioning on simplicity proves the claim. If the threshold for $r$ exceeds one, the claim is automatic.
\end{proof}

\subsection{The maximum-cut bounds at degrees fourteen and nineteen}\label{s3:app:clms}

The one-unit gap needs sharper bisection obstructions at degrees fourteen and nineteen. A single specialization of \cite[Corollary~1.3 and equations~(1.7)--(1.12)]{CLMS} gives both. Let $R_d$ be the $(d+1)\times(d+1)$ reflected-walk matrix: $(R_d)_{0,1}=(R_d)_{d,d-1}=1$, $(R_d)_{j,j-1}=(R_d)_{j,j+1}=1/2$ for $1\le j<d$, and all other entries are zero. Let $I$ be the identity matrix. For $0<\alpha<1/2$, $0<z<1$, put
\[
 B_{\alpha,z}=(1-2\alpha)I+2\alpha\sqrt z\,R_d,\quad
 \xi_j=z^{-j/2},\quad \mathbf e_0^T=(1,0,\ldots,0),
\]
where $\xi=(\xi_0,\ldots,\xi_d)^T$. Write $\operatorname{MaxCut}(G)$ for the maximum number of edges crossing a cut of $G$. The cited bound is
\[
 \frac{\operatorname{MaxCut}(G_{n,d})}{n}
 \le\frac d2+\Psi_d(\alpha,z)+o_n(1),\qquad
 \Psi_d=\frac{-\log(\mathbf e_0^TB_{\alpha,z}^d\xi)
 +(d/2)\log(1-2\alpha^2+2\alpha^2z)}{\log z}.
\]
Set $z=9/16$ and let $\alpha\uparrow1/2$. Then $B_{\alpha,z}\to(3/4)R_d$. Starting at zero, the reflected walk has position $|2J-d|$ after $d$ steps, where $J$ has distribution $\operatorname{Bin}(d,1/2)$. Hence
\[
 N_d:=\mathbf e_0^T((3/4)R_d)^d\xi
 =2^{-d}\sum_{j=0}^{d}\binom dj
       (9/16)^{\min\{j,d-j\}}.
\]
The common limiting expression is
\[
 \lim_{\alpha\uparrow1/2}\Psi_d(\alpha,9/16)
 =\frac{-\log N_d+(d/2)\log(25/32)}{\log(9/16)}.
\]
For the two required degrees this gives
\begin{align}
 \lim_{\alpha\uparrow1/2}\Psi_{14}(\alpha,9/16)
 &<-2,
 \label{eq:clms14}\\
 \lim_{\alpha\uparrow1/2}\Psi_{19}(\alpha,9/16)
 &<-3.
 \label{s3:eq:clms19}
\end{align}
Since $\log(9/16)<0$, the inequality with target $-h$ is equivalent to
\[
 N_d^2<(25/32)^d(9/16)^{2h}.
\]
The binomial sum and the rational powers satisfy
\[
\begin{array}{cc|cc}
 d&h&\text{upper bound for }N_d&\text{lower bound for }(25/32)^d(9/16)^{2h}\\\hline
 14&2&.054&.0031\\
 19&3&.0165&.00029
\end{array}
\]
The square of each entry in the third column is smaller than the entry in the fourth. Continuity gives the strict inequalities at some admissible $\alpha<1/2$. Thus a.a.s.
\[
 \operatorname{MaxCut}(G_{n,14})<(5-\epsilon)n,\qquad
 \operatorname{MaxCut}(G_{n,19})<(13/2-\epsilon)n
\]
for some fixed $\epsilon>0$. A bisection is a cut, so $p_{\rm bis}(14)\ge3$ and $p_{\rm bis}(19)\ge4$.

\subsection{Small sets}\label{s3:app:small}

The two-class count excludes small deficient sets before a signing is chosen. For $d\in\{15,18,19\}$, let $q_d$ denote the target in Lemma~\ref{s3:lem:transfer}. We prove that a.a.s.
\begin{align}
 |\cut(S)|&>(293/500)d|S|
 &&(d\ge21,\quad0<|S|\le7n/50),\label{s3:eq:smalluniform}\\
 |\cut(S)|&>q_d|S|
 &&(d\in\{15,18,19\},\quad0<|S|\le13n/100).
 \label{s3:eq:smallexception}
\end{align}
We apply Lemma~\ref{s1:lem:smallcriterion} with the three parameter triples
\[
 (d,\theta,\eta)\in
 \{(21,293/500,7/50),\ (15,267/500,13/100),\ (18,5/9,13/100)\}.
\]
For each triple, the numerator of $F_{d,\theta}''$ in~\eqref{s1:eq:smallsecond} is decreasing in $a$ and positive at $a=\eta$. Thus $F_{d,\theta}$ is convex on $[0,\eta]$. The endpoint bounds are
\[
 F_{21,293/500}(7/50),\quad
 F_{15,267/500}(13/100),\quad F_{18,5/9}(13/100)<0;
\]
the three values are approximately $-0.0070$, $-0.0223$ and $-0.0474$.
Since $F_{d,\theta}(0)=0$, the chord from zero to $\eta$ gives the negative linear bound required by the lemma. Its other parameter conditions hold for all three triples. Finally, $F_{d,\theta}$ decreases with $d$ and increases with $\theta<1-a$. The first triple therefore covers every $d\ge21$, and the second covers degree nineteen because $q_{19}/19<q_{15}/15$.

\subsection{An endpoint bound for four classes}\label{s3:app:global}

We choose weights that vary affinely on the two large classes and proportionally to size on the two small classes. The proportional weights cancel the small-class logarithms, and convexity reduces the bound to a function of the distance $t$ from a bisection side.

Fix $d$, a density threshold $c_0$, a target $q$, and $0<\eta<1/2$. Retain the four classes $(a,b,x,y)$ from Section~\ref{sec:counting}, with $x=1/2-a$ and $y=1/2-b$, and suppose
\[
 a\ge b\ge0,\qquad \eta\le a+b\le1/2,\qquad
 c_\sigma\ge c_0,\qquad |\cut(S)|/n\le q(a-b).
\]
Put
\[
 t=b+x,\quad u=x-b,\quad r_\eta=1/2-\eta,\qquad
 b_0(t)=\max\{0,(t-r_\eta)/2\}.
\]
Thus $0\le u\le\min\{t,r_\eta\}$ and $b_0(t)\le b\le t/2$. The point $b_0(t)$ is the smallest small-class size allowed by $a+b\ge\eta$. Choose
\[
 w,\kappa>0,\quad \nu\ge0,\quad U\ge1,\quad0<V\le1,
 \qquad s_0=w-\nu/2>0,
\]
and use the weights
\[
 v=(w-\nu x,\kappa b,\kappa x,w-\nu b),\qquad
 R_{ij}=v_iv_jU^{\chi_{ij}}V^{\psi_{ij}},\qquad W=\sum_{i,j}R_{ij}.
\]
The large-class weights are $s_0+\nu a$ and $s_0+\nu y$. Write $K=UV$ and
\begin{align*}
 C_-&=\frac{1-K}{2}(\nu^2+\kappa^2)+\nu\kappa(U-V),\\
 C_+&=\frac{1+K}{2}(\nu^2+\kappa^2)-\nu\kappa(U+V).
\end{align*}
Direct summation gives
\begin{equation}\label{eq:endpoint-W}
 W=2(1+K)(w-\nu t/2)^2
   +2\kappa(U+V)(w-\nu t/2)t
   +\frac{1+K}{2}\kappa^2t^2+C_-u^2.
\end{equation}
Define
\begin{align*}
 \varphi(z)&=(d-1)L(z)-dz\log(s_0+\nu z),\\
 g(t,b)&=\varphi(1/2-t+b)+\varphi(1/2-b)-dt\log\kappa
          -c_0\log U-q(1/2-t)\log V+\frac d2(W-1).
\end{align*}
Lemma~\ref{s1:lem:comparison} and the cancellation of the small-class logarithms yield
\begin{equation}\label{eq:endpoint-g}
 \mathcal F_d\le g(t,b)-L(b)-L(t-b).
\end{equation}
All formulas extend continuously when a small class is empty.

\begin{lemma}\label{lem:endpoint-reduction}
Let $0<l<h\le1/2$ and $M>1$. Put $m_0=d(s_0/w)^2-1$ and $\gamma=\log M$. Suppose
\begin{equation}\label{eq:endpoint-conditions}
 m_0>0,\qquad C_-\ge0,\qquad
 0<\gamma\le(m_0+dC_-)\min\{l,r_\eta\}.
\end{equation}
For $l\le t\le h$, define
\begin{equation}\label{eq:endpoint-R}
 R(t)=g(t,b_0(t))-L(t)+b_0(t)\log M
          +t\log(1+M^{-1}).
\end{equation}
Then $\mathcal F_d\le R(t)$. On each side of $r_\eta$, the function $R''$ is increasing. Consequently, if $R''$ is nonnegative at the left endpoint of each nonempty piece of $[l,h]\setminus\{r_\eta\}$, then
\[
 \mathcal F_d\le
 \max\{R(l),R(h),R(r_\eta)\text{ if }l<r_\eta<h\}.
\]
The one-sided derivative is used when the left endpoint is $r_\eta$.
\end{lemma}

\begin{proof}
We first bound the dependence on $b$. The assumption $m_0>0$ gives, for $0<z\le1/2$,
\[
 \varphi''(z)=\frac{ds_0^2}{z(s_0+\nu z)^2}-\frac1z\ge2m_0.
\]
The same assumption and $s_0+\nu z\le w$ imply
\[
 \varphi'''(z)=\frac{1-ds_0^2(s_0+3\nu z)/(s_0+\nu z)^3}{z^2}<0.
\]
By~\eqref{eq:endpoint-W}, for fixed $t$ the function $g(t,b)$ is symmetric about $b=t/2$ and
\[
 \frac{\partial^2g}{\partial b^2}\ge4(m_0+dC_-).
\]
Put $u_0=t-2b_0(t)=\min\{t,r_\eta\}$. Integrating the last bound from the midpoint gives
\[
 g(t,b_0(t))-g(t,t/2)\ge\tfrac12(m_0+dC_-)u_0^2.
\]
Convexity bounds $g$ by the chord joining these two values. Its slope is at most $-(m_0+dC_-)u_0\le-\gamma$, so
\[
 g(t,b)\le g(t,b_0(t))-\gamma(b-b_0(t)).
\]
The elementary entropy identity
\[
 \max_{0\le b\le t}\{-L(b)-L(t-b)-\gamma b\}
 =-L(t)+t\log(1+e^{-\gamma})
\]
follows by differentiation. Taking $\gamma=\log M$ proves the first assertion.

It remains to bound $R$ as $t$ varies. On either side of $r_\eta$, differentiation gives
\begin{equation}\label{eq:endpoint-curvature}
 R''(t)=
 \begin{cases}
 \varphi''(1/2-t)+d(\nu^2+\kappa^2-2\nu\kappa V)-1/t,
       &t<r_\eta,\\[2pt]
 \tfrac14\{\varphi''((1/2+\eta-t)/2)+
            \varphi''((3/2-\eta-t)/2)\}+dC_+-1/t,
       &t>r_\eta.
 \end{cases}
\end{equation}
Each term involving $\varphi''$ increases with $t$, because $\varphi'''<0$, and so does $-1/t$. Thus a nonnegative value at the left endpoint makes $R$ convex on that piece. Its maximum is at an endpoint, proving the lemma.
\end{proof}

\subsection{Cuts close to a bisection side}\label{s3:app:near}

Near a bisection side, the entropy gained by splitting the two small classes is dominated by a quadratic penalty. After removing this split, we compare every degree to degree twenty-one and finish with two convex functions.

Suppose $d\ge21$, $c=c_\sigma\ge d/4+(19/60)\sqrt d$, and
$|\cut(S)|/n\le(2c-2)(a-b)$. Use the comparison function $g$ in \eqref{eq:endpoint-g}, with
\[
 w=\frac12,\qquad \nu=\frac34,\qquad \kappa=\frac23,
 \qquad U=T=V^{-1}=\frac{49}{20},\qquad c_0=c,\qquad q=2c-2.
\]
Write $g_d$ and $\varphi_d$ to indicate the degree, and suppose $0\le t\le2/9$.
The curvature bound in the preceding proof gives
\[
 \partial_b^2g_d\ge4(d/16-1+dC_-)>80.
\]
Symmetry about $b=t/2$ and convexity of $g_d(t,b)-40b^2$ imply
$g_d(t,b)\le g_d(t,0)-40b(t-b)$.
For $H(u)=-u\log u-(1-u)\log(1-u)$, Jensen's inequality gives
$H(u)\le2\log(\sqrt u+\sqrt{1-u})\le2\sqrt{u(1-u)}$.
Consequently,
\begin{equation}\label{s3:eq:near-entropy}
 \mathcal F_d\le g_d(t,0)-L(t)+\frac1{40},
\end{equation}
because $2\sqrt z-40z\le1/40$ for $z\ge0$. This includes $t=0$ by continuity.

At $b=0$, put $a=1/2-t$. Formula~\eqref{eq:endpoint-W} becomes
\[
 W(t,0)=1+\{-2\nu+\kappa(T+T^{-1})\}t
           +(\nu^2+\kappa^2-2\nu\kappa/T)t^2.
\]
To compare degrees, isolate the coefficient
\[
 A(t)=a\log\frac a{1/2-\nu t}
          -t\log\kappa-\frac t2\log T+\frac12(W(t,0)-1).
\]
Then
\[
 g_d(t,0)=dA(t)-L(a)-L(1/2)
 -(1-2t)\log T-2(c-d/4)t\log T.
\]
The density assumption gives $2(c-d/4)>29/10$. Therefore, once $A\le0$ is established, the right-hand side is at most
\[
 21A(t)-L(a)-L(1/2)-(1+9t/10)\log T.
\]
For $\theta=1/6$, the tangent bound $-L(t)\le\theta-t(1+\log\theta)$ now reduces~\eqref{s3:eq:near-entropy} to $\mathcal F_d\le B(t)$, where
\[
\begin{split}
 B(t)={}&21A(t)-L(a)-L(1/2)-(1+9t/10)\log T\\
        &+\theta-t(1+\log\theta)+\frac1{40}.
\end{split}
\]

Both required bounds follow from convexity. Indeed,
\[
 B''(t)=\varphi_{21}''(a)+21(\nu^2+\kappa^2-2\nu\kappa/T)>0,
 \qquad 21A''(t)=B''(t)+1/a.
\]
The endpoint inequalities
\[
 A(0)=0,\qquad A(2/9)<-\frac1{4000},\qquad
 B(0),B(2/9)<-\frac1{100}
\]
give $A\le0$ and $B<-1/100$ on $[0,2/9]$. Hence
\begin{equation}\label{s3:eq:near}
 \mathcal F_d<-1/100.
\end{equation}

\subsection{Three comparisons for all degrees}\label{s3:app:tables}\label{app:fixedtargets}

Two choices of weights cover the remaining sets at degrees fifteen, eighteen and nineteen. A third covers every $d\ge21$. In each case we compare degrees through the coefficient of $d$ in the exponent.

\begin{table}[htbp]
\centering\small
\setlength{\tabcolsep}{5pt}
\caption{The three choices of comparison weights. Every decimal is exact.}\label{tab:endpoint-parameters}
\begin{tabular}{clrrrrrr}
\toprule
Row&$t$-range&$w$&$\nu$&$\kappa$&$U$&$V$&$M$\\
\midrule
I&$[3/25,3/10]$&.43&.45&.63&3.10&.27&8\\
II&$[3/10,1/2]$&.42&.44&.73&2.50&.24&32\\
III&$[2/9,1/2]$&.425&.64&.78&2.12&$U^{-1}$&32\\
\bottomrule
\end{tabular}
\end{table}

\paragraph{The finite degrees.}
Put $\eta=13/100$ and $r=1/2-\eta=37/100$. Let $c_0,q$ be the targets in Lemma~\ref{s3:lem:transfer}. The distance inequality in Section~\ref{sec:counting} gives
\[
 t\ge l_d:=\frac{c_0-q/2}{d-q},\qquad
 l_{15}=\frac{197}{1398},\quad l_{18}=\frac18,\quad
 l_{19}=\frac{54}{445}.
\]
All three exceed $3/25$. Write $R_{d,c,q}$ for the function in~\eqref{eq:endpoint-R}, with weights from the row under consideration.

We first obtain $\mathcal F_d\le R_{d,c_0,q}$ from the first assertion of Lemma~\ref{lem:endpoint-reduction}, using row I on $[l_d,3/10]$ and row II on $[3/10,1/2]$. At degree fifteen both rows have $m_0>12/5$ and $C_->17/20$, and $m_0$ increases with $d$. Together with $l_{15}>7/50$ and $l_{18},l_{19}>3/25$, these bounds imply the slope condition~\eqref{eq:endpoint-conditions} for every actual triple $(d,c_0,q)$.

We bound these functions together. Set
\[
 \bar c_d=\frac{499}{100}+\frac{11}{40}(d-15),\qquad
 \bar q_d=\frac{801}{100}+\frac{21}{40}(d-15),\qquad
 R_{15}=R_{15,499/100,801/100}.
\]
The affine dependence on the parameters gives
\begin{equation}\label{eq:shared-affine}
 R_{d,\bar c_d,\bar q_d}(t)=R_{15}(t)+(d-15)A(t),
\end{equation}
where, with $b=b_0(t)$, $x=t-b$, $a=1/2-x$ and $y=1/2-b$,
\[
\begin{split}
 A(t)={}&L(a)+L(y)-a\log(w-\nu x)-y\log(w-\nu b)
       -t\log\kappa+\tfrac12(W-1)\\
       &-\tfrac{11}{40}\log U-\tfrac{21}{40}(1/2-t)\log V.
\end{split}
\]
For the actual targets,
\begin{equation}\label{eq:shared-domination}
 R_{d,c_0,q}(t)\le R_{d,\bar c_d,\bar q_d}(t).
\end{equation}
There is equality at degree fifteen. At degree nineteen, $c_0=\bar c_{19}+1/25$ and $q=\bar q_{19}-1/100$, so both changes decrease $R$. At degree eighteen the difference is
\[
 -\frac{37}{200}\log U-\frac{83}{200}(1/2-t)\log V<0.
\]
For row I this follows from $U^7V^6>1$ and $1/2-t\le19/50$; for row II it follows from $U^2V>1$ and $1/2-t\le1/5$.

It remains to bound the two degree-fifteen functions $R_{15}$ and $A$. The function $R_{15}$ bounds $\mathcal F_{15}$ only for $t\ge l_{15}$. We nevertheless bound it on the whole interval $[3/25,3/10]$ of row~I, because \eqref{eq:shared-affine} expresses the functions at degrees eighteen and nineteen, whose domains begin below $l_{15}$, through $R_{15}$ and $A$. Formula~\eqref{eq:endpoint-curvature} gives $R_{15}''>7$ at $3/25$ for row I and at $3/10,r+$ for row II. Since $R_{15}''$ increases on each side of $r$, it is positive throughout each piece. Moreover,
\[
 15A''=R_{15}''+\frac{d^2}{dt^2}\{L(a)+L(y)+L(t)\}>0.
\]
Hence both functions are convex on each piece. Their endpoint bounds are
\begin{equation}\label{eq:endpoint-negative}
 R_{15}(t),\ A(t)<-\frac1{200}
 \quad\begin{cases}
 t\in\{3/25,3/10\},&\text{row I},\\
 t\in\{3/10,37/100,1/2\},&\text{row II}.
 \end{cases}
\end{equation}
Convexity extends these bounds to the full intervals. Equations~\eqref{eq:shared-affine}--\eqref{eq:shared-domination} now give $\mathcal F_d<-1/200$ throughout the remaining region at all three finite degrees.

\paragraph{The uniform range.}
Use row III with
\[
 d=21,\qquad \eta=7/50,\qquad c_0=67/10,\qquad q=57/5=2c_0-2.
\]
Appendix~\ref{s3:app:near} already excludes $t\le2/9$. For the remaining interval, $m_0>1/6$, $C_->0$, and the slope condition~\eqref{eq:endpoint-conditions} holds. The curvature is positive at $2/9$ and $9/25+$, while
\[
 R(2/9),\ R(9/25),\ R(1/2)<-1/1000.
\]
Lemma~\ref{lem:endpoint-reduction} gives the same bound throughout the interval.

We now pass from degree $21$ to all larger degrees. For row III put $T=U=V^{-1}$, and define
\begin{equation}\label{s3:eq:A}
 \mathscr A=L(a)+L(y)-a\log(w-\nu x)-y\log(w-\nu b)
             -t\log\kappa-\tfrac12t\log T+\tfrac12(W-1).
\end{equation}
If the actual bisection density is $c$, the two cut terms in \eqref{s1:eq:comparison} are at most
\[
 -c\log T+(2c-2)(1/2-t)\log T
 =-2ct\log T-(1-2t)\log T.
\]
This decreases with $c$. Substituting $c\ge d/4+(19/60)\sqrt d$ therefore gives
\begin{equation}\label{s3:eq:G}
 \mathcal F_d\le d\mathscr A-L(a)-L(y)-L(b)-L(x)
       -(1-2t)\log T-\frac{19}{30}\sqrt d\,t\log T.
\end{equation}
It suffices to prove $\mathscr A<0$, for then the right-hand side decreases with $d$. The function $z\mapsto L(z)-z\log(s_0+\nu z)$ has positive second derivative. Together with $C_-\ge0$, this makes $\mathscr A$ convex and symmetric about $b=t/2$ for fixed $t$. Its maximum on $[b_0(t),t/2]$ is therefore at $b_0(t)$.

Write $A_0(t)=\mathscr A(t,b_0(t))$ and
$J_0(t)=L(1/2-t+b_0(t))+L(1/2-b_0(t))$. On each piece, the degree-$21$ expression~\eqref{eq:endpoint-R} satisfies
\[
 21A_0''(t)=R''(t)+J_0''(t)+1/t>0.
\]
Thus $A_0$ is convex there. At its three endpoints, substitution gives
\[
 A_0(2/9),\ A_0(9/25),\ A_0(1/2)<-1/20000.
\]
Consequently $\mathscr A<0$ throughout the domain. The right-hand side of~\eqref{s3:eq:G} decreases with $d$, and at $d=21$ it is at most $g(t,b)-L(b)-L(x)$ for row III, since $c_0=67/10<21/4+(19/60)\sqrt{21}$. Lemma~\ref{lem:endpoint-reduction} bounds this value by $-1/1000$. This proves the uniform estimate for every $d\ge21$.

The endpoint values, rounded to four decimal places, are listed below. Appendix~\ref{app:arithmetic} gives the logarithm bounds used to verify the rational estimates above.
\begin{center}\small
\begin{tabular}{@{}cccc@{\qquad}cccc@{}}\toprule
row&$t$&$R_{15}(t)$&$A(t)$&row&$t$&$R(t)$&$A_0(t)$\\\midrule
I&$3/25$&$-0.0078$&$-0.0053$&III&$2/9$&$-0.0091$&$-0.0072$\\
I&$3/10$&$-0.0340$&$-0.0104$&III&$9/25$&$-0.0047$&$-0.0005$\\
II&$3/10$&$-0.0475$&$-0.0222$&III&$1/2$&$-0.4878$&$-0.0263$\\
II&$37/100$&$-0.0287$&$-0.0168$&&&&\\
II&$1/2$&$-0.4712$&$-0.0537$&&&&\\\bottomrule
\end{tabular}
\end{center}

Together with the near-side and small-set estimates, this completes the proof of Lemma~\ref{s3:lem:transfer}.

\section{The quantitative bisection estimate}\label{app:dms}

We prove Lemma~\ref{s2:lem:DMS} by transferring the Gaussian comparison from a graph with independently sampled edges to the regular pairing model. We first pass to independent Poisson degrees and truncate degrees above $d$. We then restore regularity, subtracting the mean increase in cut size and bounding the fluctuations uniformly over balanced cuts. The Gaussian comparison contributes the leading error, of order $d^{1/3}$.

\subsection{Independent edges and Poisson degrees}

For even $n$, let $\Omega_n$ be the balanced signings of $[n]$, and define
\[
 a_n=\frac1n\E\max_{\tau\in\Omega_n}
       \frac1{\sqrt n}\sum_{1\le i<j\le n}g_{ij}\tau_i\tau_j,
\]
where the $g_{ij}$ are independent Gaussian variables of mean zero and variance one. This agrees with the expected-maximum normalization in \cite[equation~(1.2)]{DMS}: the two Gaussian coefficients of each unordered pair combine to give variance $1/n$, while the diagonal term has mean zero and is constant over the signing. The balancing argument in \cite[Lemma~2.4 and equations~(2.12)--(2.13)]{DMS} gives $a_n\to C^*$.

Fix $\nu\ge2$. Let $\mathcal P_{n,\nu}$ have a Poisson number of edges of mean $\nu n/2$, each with two independently uniform endpoints in $[n]$. Loops and repeated edges are allowed. Cut sizes count edges with multiplicity, and loops never cross a cut. The dilute model of \cite{DMS} has mean edge count $\gamma n$, so its parameter here is $\gamma=\nu/2$. We use the second inequality in \cite[Lemma~2.3, equation~(2.10)]{DMS}, for the maximum energy constrained to balanced signings. In its notation, $H^D_\gamma(\sigma)=-\sum_{e=ij}\sigma_i\sigma_j$, $\widehat e^D_n(\gamma)=n^{-1}\E\max_{\sigma\in\Omega_n}H^D_\gamma(\sigma)$, and the constrained Gaussian expected minimum energy per vertex satisfies $e^{\rm SK}_n=-a_n$. A cut has size $(|E|+H^D_\gamma(\sigma))/2$. Taking the absolute constant $A_0$ from that lemma with $\gamma_0=1$, we obtain
\[
 \E c(\mathcal P_{n,\nu})
 =\frac\nu4+\frac12\widehat e^D_n(\nu/2),\qquad
 \left|\frac{\widehat e^D_n(\nu/2)}{\sqrt\nu}-a_n\right|
 \le A_0(\nu/2)^{-1/6}.
\]
It follows that, for an absolute constant $K_{\rm int}$,
\begin{equation}\label{app:eq:diluteexpectation}
 \left|\E c(\mathcal P_{n,\nu})-\frac\nu4-\frac{\sqrt \nu}{2}a_n\right|
 \le K_{\rm int}\nu^{1/3}.
\end{equation}
The concentration estimate \cite[Lemma~2.1]{DMS} now yields, at every fixed $\nu$,
\begin{equation}\label{app:eq:dilutebound}
 \left|c(\mathcal P_{n,\nu})-\frac\nu4-\frac{C^*}{2}\sqrt \nu\right|
 \le K_{\rm int}\nu^{1/3}+o_n(1)
 \qquad\text{a.a.s.}
\end{equation}
The truncation argument is easier with independent degrees. Let $\mathcal H_{n,\nu}$ be formed by giving vertex $i$ an independent Poisson number $X_i$ of labelled points, of mean $\nu$, and matching all points uniformly. Whenever a point total is odd, leave one uniformly chosen point unmatched; use the same convention in every rematching below. Uniformity in the odd case follows by permutation symmetry. Changing the unmatched point can be coupled by changing at most two edges and thus contributes only $O(1/n)$ to a cut density.

To couple the two models, take an infinite sequence of independent uniform vertices of $[n]$. Independently draw Poisson variables $N,M$ of means $\nu n$ and $\nu n/2$. The first $N$ entries give independent Poisson vertex counts of mean $\nu$. Pairing consecutive entries, and ignoring the last entry when $N$ is odd, has the conditional law of a uniform pairing with those counts, so it generates $\mathcal H_{n,\nu}$. The first $2M$ entries generate $\mathcal P_{n,\nu}$. The graphs differ by $|\lfloor N/2\rfloor-M|$ edges, whose expectation is $O_\nu(\sqrt n)$. This is $o(n)$ a.a.s. at fixed $\nu$. Changing one edge changes the maximum bisection by at most one, and hence
\begin{equation}\label{app:eq:poissonpairing}
 |c(\mathcal H_{n,\nu})-c(\mathcal P_{n,\nu})|=o_n(1)
 \qquad\text{a.a.s.}
\end{equation}

\subsection{Truncating the degrees}

Fix a sufficiently large integer $d$, and put
\[
 \Delta=\sqrt d\log d,\qquad \nu=d-\Delta\ge2.
\]
Let $X$ be Poisson of mean $\nu$, write $(x)_+=\max\{x,0\}$, and put $\delta_d=\E(X-d)_+$. Retain $\min\{X_i,d\}$ of the points at vertex $i$. Delete edges meeting discarded points and uniformly rematch the retained points left free. Write $G_2$ for the resulting graph. Its pairing is uniform conditional on the retained degrees, by \cite[Lemma~3.4 and its application in Section~3.1]{DMS}, with the unmatched-point convention above when necessary.

If $K=\sum_i(X_i-d)_+$ points are discarded, at most $K$ old edges are removed and at most $(K+1)/2$ new edges are added. Thus the difference of the two cut sizes, for every vertex partition, is at most $3K/2+1$. The law of large numbers gives
\begin{equation}\label{app:eq:truncation}
 |c(G_2)-c(\mathcal H_{n,\nu})|
 \le\tfrac32\delta_d+o_n(1)
 \qquad\text{a.a.s.}
\end{equation}
Since $\delta_d\le\E[X_i\mathbf 1_{\{X_i\ge d+1\}}]=\nu\Pr(X_i\ge d)$, the Poisson exponential-moment bound gives
\begin{equation}\label{app:eq:poissontail}
 \delta_d\le d\exp\{-d\log(d/\nu)+d-\nu\}
 \le d\exp\{-\tfrac12(\log d)^2\}
 =:\omega_d=o_d(1).
\end{equation}
The second inequality uses $-\log(1-u)-u\ge u^2/2$ at $u=\Delta/d$.

\subsection{Recovering the regular graph}

When restoring the missing degree $\Delta=\sqrt d\log d$, we subtract the leading mean contribution $\Delta/4$ to the cut density and bound the remaining error uniformly over balanced signings.

Start with $d$ points at every vertex. Using Poisson counts $X_i$ as above, colour
\[
 Z_i=(d-X_i)_+
\]
of its points blue, leaving the other $\min\{X_i,d\}$ points red. Choose the point labels of each colour uniformly, and pair all $dn$ points uniformly, independently of the colours. The resulting graph $G_1$ is the $d$-regular pairing model.

Let $F_{\mathrm D}$ be the family of edges meeting a blue point, and let $F_{\mathrm{BB}}\subseteq F_{\mathrm D}$ consist of the blue--blue edges. Delete $F_{\mathrm D}$ and uniformly rematch the freed red points; write $F_{\mathrm Q}$ for the new edges. Together with the unchanged red--red edges, these form a graph with the conditional law of $G_2$ from the preceding subsection. Thus the two constructions of $G_2$ may be coupled.

For an edge family $F$, let $Q_F(\sigma)$ count its edges crossing the bisection of a balanced signing $\sigma$. Then
\begin{equation}\label{app:eq:rewiringidentity}
 Q_{G_1}(\sigma)-Q_{G_2}(\sigma)
 =Q_{F_{\mathrm D}}(\sigma)-Q_{F_{\mathrm Q}}(\sigma).
\end{equation}

We first compute the mean of this difference for a fixed signing. Put $\mu=\E Z_i=\Delta+\delta_d$ and $B=\sum_iZ_i$. Counting blue points gives
\[
 |F_{\mathrm D}|=B-|F_{\mathrm{BB}}|,
 \qquad |F_{\mathrm Q}|=\left\lfloor\frac{B-2|F_{\mathrm{BB}}|}{2}\right\rfloor.
\]
The inequality $|(d-x)_+-(d-y)_+|\le|x-y|$ gives $\operatorname{Var}(Z_i)\le\operatorname{Var}(X_i)=\nu\le d$. Each pair of points is matched with probability $1/(dn-1)$, and therefore
\[
 \E|F_{\mathrm{BB}}|=\frac{\E[B(B-1)]}{2(dn-1)}
 =O(n\mu^2/d)+O_d(1).
\]
The deleted and rematched families are invariant under vertex relabelling and each has $O_d(1)$ expected loops. Every non-loop edge crosses a fixed balanced bisection with probability $n/[2(n-1)]$, so
\[
 \frac1n\E Q_{F_{\mathrm D}}(\sigma)
 =\frac\mu2+O(\mu^2/d)+o_n(1),\qquad
 \frac1n\E Q_{F_{\mathrm Q}}(\sigma)
 =\frac\mu4+O(\mu^2/d)+o_n(1).
\]
Here $\mu^2/d=O((\log d)^2)$, and replacing $\mu$ by $\Delta$ costs $O(\delta_d)$.

To allow a signing chosen from the graph, apply the uniform fluctuation estimate \cite[Lemma~3.3, equation~(3.7)]{DMS} with $\gamma=d$ and $\gamma_-=\nu$. It gives an absolute constant $K_{\rm fl}$ such that, for $F=F_{\mathrm D}$ or $F_{\mathrm Q}$,
\[
 \Pr\left(
 \sup_{\sigma\in\Omega_n}
 \bigl|Q_F(\sigma)-\E Q_F(\sigma)\bigr|
 >K_{\rm fl}nd^{1/4}\sqrt{\log d}
 \right)\longrightarrow0.
\]
The expectations average both colours and matchings. Combining this estimate with the two means and \eqref{app:eq:rewiringidentity} gives an absolute constant $K_{\rm rw}$ such that
\begin{equation}\label{app:eq:uniformrewiring}
 \sup_{\sigma\in\Omega_n}
 \left|\frac{Q_{G_1}(\sigma)-Q_{G_2}(\sigma)}n-\frac\Delta4\right|
 \le K_{\rm rw}\{d^{1/4}\sqrt{\log d}+(\log d)^2+\delta_d\}+o_n(1)
 \qquad\text{a.a.s.}
\end{equation}

\subsection{Combining the estimates}

The simultaneous bound~\eqref{app:eq:uniformrewiring} remains valid when the two cut counts are replaced by their separate maxima. Combine it with the independent-edge estimate~\eqref{app:eq:dilutebound}, the coupling~\eqref{app:eq:poissonpairing}, and the truncation bound~\eqref{app:eq:truncation}. The mean terms add to $(\nu+\Delta)/4=d/4$, while replacing $\sqrt\nu$ by $\sqrt d$ costs $O(\log d)$. Hence, for an absolute constant $K_{\rm reg}$,
\[
 \left|c(G_1)-\frac d4-\frac{C^*}{2}\sqrt d\right|
 \le K_{\rm int}d^{1/3}
 +K_{\rm reg}\{d^{1/4}\sqrt{\log d}+(\log d)^2+\omega_d+\log d\}
 +o_n(1).
\]
The terms other than $K_{\rm int}d^{1/3}$ are $o_d(d^{1/3})$. Absorbing them into the constant proves Lemma~\ref{s2:lem:DMS} for the regular pairing model. Conditioning on simplicity, whose probability has a positive limit at each fixed $d$, gives the conclusion for $G_{n,d}$.

\section{Directed cuts in marked orientation models}\label{app:marked-expansion}

We prove the directed cut bounds in Proposition~\ref{prop:small-strong} by a matching count. At degree nine we retain the depletion of available points, while a coarser estimate suffices at degrees thirteen and seventeen. The two subsections then treat Hamilton augmentation at degree eleven and matching augmentation at degree seven.

\begin{lemma}[Expansion in the marked orientation model]\label{lem:marked-expansion}
Let $p\in\{2,3,4\}$ and $d=4p+1$. Put $h=1/40,1/20,1/10$ for $p=2,3,4$, respectively. Give half the vertices indegree $p$ and outdegree $3p+1$, and the other half the opposite degrees. Match all outgoing points uniformly to all incoming points. Asymptotically almost surely, every nontrivial vertex set $S$ satisfies
\[
 \min\{d^+(S),d^-(S)\}>h\min\{|S|,n-|S|\},
\]
where $d^+(S)$ and $d^-(S)$ count the arcs leaving and entering $S$.
\end{lemma}
\begin{proof}
For the counting step, allow any $0<h<p$. By complementation and reversal, it suffices to count sets $S$ with $|S|=sn\le n/2$ and $d^+(S)\le hsn$. Suppose $S$ contains $an$ vertices of outdegree $p$ and $bn$ of outdegree $3p+1$, where $a+b=s$, and put $j=\lfloor hsn\rfloor$.

We may assume $a\ge b$, at the cost of a factor two in the union bound. Indeed, the failure probability is nonincreasing in the outgoing volume and nondecreasing in the incoming volume of $S$, by coupling uniform samples of the incoming points. Reflecting a pair with $a<b$ preserves the number of sets and increases its failure probability.

Put
\[
 t=\frac{a-b}{s},\qquad \kappa=\frac{2p+1}{d},\qquad
 r=\frac{d-(2p+1)t}{2},\qquad k=d-r,\qquad m=\frac d2.
\]
Then $0\le t\le1$. The set has $rsn$ outgoing and $ksn$ incoming points, out of $mn$ of each kind in total. Write
$\beta=ks/m=s(1+\kappa t)$ and $(z)_i=z(z-1)\cdots(z-i+1)$.
If at most $j$ arcs leave $S$, some $rsn-j$ outgoing points are all matched internally. Consequently
\begin{equation}\label{eq:marked-small-boundary}
\begin{split}
 \Pr(d^+(S)\le hsn)
 &\le\binom{rsn}{j}\frac{(ksn)_{rsn-j}}{(mn)_{rsn-j}}
 \le\binom{rsn}{j}\beta^{rsn-j}\\
 &\le\exp\{sn[h\log(er/h)+(r-h)\log\beta]\}.
\end{split}
\end{equation}
The last bound follows from the binomial estimate and is valid also when $j=0$. We will use its coarse form at $p=3,4$ and its falling-factorial form at $p=2$.

Let $H(x)=-x\log x-(1-x)\log(1-x)$, with continuous endpoint values. The number of sets of the given type is at most
\[
 \exp\left\{\frac n2[H(s(1+t))+H(s(1-t))]\right\}.
\]
Both estimates below use
\[
 H(x)\le-x\log x+x-\frac{x^2}{2}-\frac{x^3}{6}
 \qquad(0\le x\le1).
\]
Put $g(t)=((1+t)\log(1+t)+(1-t)\log(1-t))/2$, with $0\log0=0$.

\paragraph{The coarse estimate.}
For $p=3,4$, the exponent per vertex obtained from~\eqref{eq:marked-small-boundary} is
\[
 F_{p,h}(s,t)=\tfrac12[H(s(1+t))+H(s(1-t))]
       +s[h\log(er/h)+(r-h)\log\beta].
\]
The cubic entropy bound gives
\begin{align*}
 \frac{F_{p,h}(s,t)}s
 \le{}&(r-1-h)\log s+1-g(t)+(r-h)\log(1+\kappa t)
       +h\log(er/h)\\
 &-\frac{s(1+t^2)}2-\frac{s^2(1+3t^2)}6.
\end{align*}
For the stated choices of $h$, this upper bound increases with $s$ on $(0,1/2]$. Indeed, its derivative is at least $1/3-2h>0$, using $r\ge2$. Setting $s=1/2$ therefore gives
\[
 \frac{F_{p,h}(s,t)}s
 \le U_p(t)+h\log\frac{2er}{h(1+\kappa t)},
\]
where
\[
 U_p(t)=\frac{17}{24}-\frac{3t^2}{8}-g(t)
        -(r-1)\log2+r\log(1+\kappa t).
\]
Also
\[
 U_p''(t)=-\frac34-\frac1{1-t^2}
          -\frac{(2p+1)\kappa}{1+\kappa t}
          -\frac{r\kappa^2}{(1+\kappa t)^2}
 \le-\frac74\qquad(0\le t<1).
\]
In both cases the additional term is at most $h\log(17e/h)$, since $2r\le d\le17$ and $1+\kappa t\ge1$.

\paragraph{Degree nine.}
Now let $p=2$ and $h=1/40$, and put $\rho=r-h$. The falling-factorial bound in~\eqref{eq:marked-small-boundary}, with $m=9/2$ and $k=9-r$, gives the expected-count exponent
\begin{align*}
 \mathcal E(s,t)={}&
 \tfrac12\{H(s(1+t))+H(s(1-t))\}+srH(h/r)\\
 &+L(ks)-L((k-\rho)s)+L(m-\rho s)-L(m).
\end{align*}
More precisely, the logarithm of the expected count is at most $n\mathcal E(s,t)+O(\log n)$, uniformly over the class sizes. This includes the error from replacing $j$ by $hsn$.

Again we reduce to $s=1/2$. Put
\[
 T(x)=(1-x)\log(1-x)+x
     =\sum_{j\ge2}\frac{x^j}{j(j-1)}.
\]
The cubic entropy bound makes $\mathcal E(s,t)/s$ at most
\begin{align*}
 &(\rho-1)\log s+1-g(t)+rH(h/r)
       +L(k)-L(k-\rho)-\rho(\log m+1)\\
 &\hspace{15mm}-\frac{s(1+t^2)}2
       -\frac{s^2(1+3t^2)}6+\frac msT(\rho s/m).
\end{align*}
The last term increases with $s$, as its power series has positive coefficients. The derivative of the preceding terms is at least $1/3-2h>0$, since $\rho\ge2-h$. Evaluating at $s=1/2$ therefore gives $\mathcal E(s,t)/s\le V_h(t)$, where
\begin{align*}
 V_h(t)={}&\frac{17}{24}+\log2-\frac{3t^2}{8}-g(t)
       +L(k)+L(k+h)-L(5t+h)-L(9)\\
       &+L(r)-L(h)-L(r-h).
\end{align*}
As in the coarse estimate, this function is uniformly concave:
\[
 V_h''(t)=-\frac34-\frac1{1-t^2}
 +\frac{25}{4}\left\{
 \frac1k+\frac1{k+h}-\frac4{5t+h}
 +\frac1r-\frac1{r-h}\right\}\le-\frac74.
\]
Here $1/k+1/(k+h)<4/(5t+h)$, and the last two terms in braces have negative sum.

\paragraph{Completing the count.}
If $\Psi''\le-7/4$, Taylor's inequality at any $\tau\in(0,1)$ gives
\[
 \sup_{0\le t\le1}\Psi(t)
 \le\Psi(\tau)+\frac27\Psi'(\tau)^2.
\]
We apply this to $\Psi=V_{1/40},U_3,U_4$ at the respective degrees. Let $\Delta_2=0$ and $\Delta_p=h\log(17e/h)$ for $p=3,4$, with their respective values of $h$. The required scalar bounds are
\[
\begin{array}{c|c|c|c|c|c}
 p&\Psi&\tau&\Psi(\tau)<&|\Psi'(\tau)|<&\Delta_p\le\\\hline
 2&V_{1/40}&7/8&-1/80&1/15&0\\
 3&U_3&91/100&-39/100&1/50&7/20\\
 4&U_4&15/16&-7/10&1/2&31/50
\end{array}
\]
The concavity bound, followed by addition of $\Delta_p$, yields
\begin{equation}\label{eq:marked-exponent-margin}
 \mathcal E(s,t)<-s/100\quad(p=2),\qquad
 F_{p,h}(s,t)<-s/200\quad(p=3,4).
\end{equation}

There are $O(n^2)$ pairs of class sizes. These bounds, including the uniform $O(\log n)$ error for $p=2$, exclude all sets with $|S|\ge C\log n$ a.a.s., for a sufficiently large absolute constant $C$. For smaller sets, the incoming volume is at most $(3p+1)|S|$ and the outgoing volume is at least $p|S|$. Applying the coarse bound in~\eqref{eq:marked-small-boundary} and summing over the vertices gives, for a constant $C_p$,
\[
 \Pr\bigl(\exists S:\ 1\le |S|\le C\log n,\ d^+(S)\le h|S|\bigr)
 \le\sum_{1\le i\le C\log n}
       [C_p(i/n)^{p-h-1}]^i=o(1),
\]
since $p-h-1>0$. Reversal and complementation complete the proof.
\end{proof}

\subsection{Hamilton augmentation at degree eleven}\label{app:eleven-entropy}

\begin{lemma}\label{lem:eleven-marked}
Let $(d,q,\beta)=(9,5,101/200)$. Choose uniformly a degree-$d$ pairing together with one of its Jaeger orientations, of divergence $q\sigma$, and adjoin an independent uniform Hamilton cycle on the same vertices with a uniformly chosen cyclic direction. Then a.a.s. the union satisfies~\eqref{eq:eleven-cut}.
\end{lemma}

\begin{proof}
Write $P$ for the marked pairing and $C$ for the directed cycle. Fix the vertex types. Complementation gives sets of size $sn\le n/2$. Reversing all arcs and interchanging the two equally large type classes preserves the law, so it suffices, up to a factor two, to count sets containing $(1+t)sn/2$ vertices of outdegree $(d-q)/2$ and $(1-t)sn/2$ of outdegree $(d+q)/2$, where $0\le t\le1$. Put
\[
 A=(d-qt)/2,\qquad B=(d+qt)/2,\qquad M=d/2.
\]
The outgoing and incoming point totals of the set in $P$ are $Asn$ and $Bsn$. If $J$ and $K$ count its outgoing boundary arcs in $P$ and $C$, then
\[
 |\cut(S)|=qtsn+2(J+K).
\]
Failure of~\eqref{eq:eleven-cut} therefore requires $J+K\le\beta tsn$.

For a fixed set of $m$ vertices, the probability that $K=k>0$ is
\[
 \frac{(n/k)\binom{m-1}{k-1}\binom{n-m-1}{k-1}}{\binom nm}.
\]
Indeed, its membership indicators around the cycle have $k$ runs of each type. The matching probability is hypergeometric. Combining these two counts with the choices of the two vertex classes gives the exponent $F(s,t,r)$ defined below, with $J=rsn$ and $K=(\beta t-r)sn$.

It is enough to consider this equality boundary. The two rates are concave in their separate boundary sizes and are maximized at their asymptotic means. Divided by $sn$, the sum of those means is
\[
 A(1-Bs/M)+1-s\ge A^2/d+1/2>\beta t.
\]
Thus a maximum under $J+K\le\beta tsn$ lies on the boundary. Lemma~\ref{lem:eleven-entropy} gives $F(s,t,r)<-\eta s$ for a fixed $\eta>0$. The factors from Stirling's error and from summing over class-size and boundary patterns are polynomial in $n$, so this excludes all deficient sets with at least $C_1\log n$ vertices, for a sufficiently large constant $C_1$.

For smaller sets of size $k$, there are at least $3k$ outgoing points in the factor and cycle together. Deficiency therefore forces at least $(3-\beta)k$ of these arcs to stay inside the set. There are only $C_2^k$ ways to select their outgoing points. The probability that $\ell$ selected matching arcs all end inside the set is at most $(C_3k/n)^\ell$. For $\ell$ selected cycle arcs, sum over at most $k^\ell$ assignments of their heads inside the set. A prescription with distinct heads and no directed cycle has probability $(n-\ell-1)!/(n-1)!$. All other prescriptions are impossible on a proper vertex set. This gives the same bound. Independence, followed by a union bound over the vertex sets, gives
\[
 \Pr\bigl(\exists S:\ 1\le |S|\le C_1\log n,
                   J+K\le\beta|S|\bigr)
 \le\sum_{1\le k\le C_1\log n}
       [C_4(k/n)^{2-\beta}]^k=o(1).
\]
Thus the marked augmented model satisfies~\eqref{eq:eleven-cut} a.a.s.
\end{proof}

It remains to prove the entropy estimate used in this proof. Take $(d,q,\beta)=(9,5,101/200)$ and put $M=d/2$. For $0<s\le1/2$, $0\le t\le1$ and $0\le r\le\beta t$, write
\[
 A=(d-qt)/2,\qquad B=(d+qt)/2,\qquad y=\beta t-r.
\]
Here $sn$ is the set size, $tsn$ its signing imbalance, and $rsn,ysn$ its outgoing boundaries in the original factor and cycle. Retain $L$ and $H$ from the preceding count. For a fixed set, the probability that the original outgoing boundary equals $j$ is
\[
 \frac{\binom{Bsn}{Asn-j}\binom{Mn-Bsn}{j}}{\binom{Mn}{Asn}}.
\]
Combining this with the cycle count and the choices of the two vertex classes gives the exponent
\begin{align*}
 F(s,t,r)={}&\tfrac12\{H(s(1+t))+H(s(1-t))\}\\
 &+L(Bs)+L(M-Bs)+L(As)+L(M-As)\\
 &-L((A-r)s)-L(rs)-L((qt+r)s)-L(M-(B+r)s)-L(M)\\
 &+2L(s)+2L(1-s)-L((1-y)s)-2L(ys)-L(1-(1+y)s).
\end{align*}

\begin{lemma}\label{lem:eleven-entropy}
Throughout this domain, $F(s,t,r)<-3s/200$.
\end{lemma}
\begin{proof}
Put $f=F/s$.

\paragraph{Curvature in the imbalance and boundary.}
The Hessian $\nabla^2 f$ below is with respect to $(t,r)$. Matrix inequalities refer to quadratic forms; for a column vector $z$, write $z^2=zz^T$. Set
\begin{align*}
 T&=\frac{q^2}{4}\left(\frac1A+\frac1B+\frac{s}{M-As}+\frac{s}{M-Bs}\right),\\
 P&=\frac12\left(\frac1{1+t}+\frac1{1-t}
       +\frac{s}{1-(1+t)s}+\frac{s}{1-(1-t)s}\right),\\
 a_1&=\frac1{A-r}+\frac{s}{M-(B+r)s},\qquad
 b_1=\frac1{1-y}+\frac2y+\frac{s}{1-(1+y)s}.
\end{align*}
Differentiation in the interior gives
\[
 \nabla^2f=(T-P)(1,0)^2-a_1(q/2,1)^2
 -r^{-1}(0,1)^2-(qt+r)^{-1}(q,1)^2-b_1(\beta,-1)^2.
\]
The positive term has the cancellation
\[
 T-P=h(t)+\frac{s}{1-s}h\left(\frac{st}{1-s}\right),\qquad
 h(z)=\frac{q^2/d}{1-q^2z^2/d^2}-\frac1{1-z^2}.
\]
Since $q^4<d^3$, differentiation shows that $h$ decreases on $[0,1)$.
Thus $T-P\le2(q^2/d-1)$. Also
\[
 a_1\ge2/d,\quad r^{-1}\ge\beta^{-1},\quad
 (qt+r)^{-1}\ge(q+\beta)^{-1},\quad
 b_1\ge(1-\beta)^{-1}+2/\beta.
\]
The last bound uses the decrease of $(1-y)^{-1}+2/y$ for $0<y\le\beta<2-\sqrt2$. It follows that $-\nabla^2f\ge K$, where
\[
 K=-2(q^2/d-1)(1,0)^2+\frac2d(q/2,1)^2+\beta^{-1}(0,1)^2
 +(q+\beta)^{-1}(q,1)^2+
 \bigl((1-\beta)^{-1}+2/\beta\bigr)(\beta,-1)^2.
\]
The diagonal entries of $K$ exceed $19/5$ and $8$, and its off-diagonal
entry has absolute value less than $8/5$. Since $(19/5-3)(8-3)>(8/5)^2$,
the principal-minor criterion gives
\begin{equation}\label{eq:eleven-hessian}
 \nabla^2f\le-3I.
\end{equation}

\paragraph{Reduction to one variable.}
At $(t_0,r_0)=(9/10,17/100)$, put $g(s)=f(s,t_0,r_0)$,
$u(s)=f_t(s,t_0,r_0)$ and $v(s)=f_r(s,t_0,r_0)$. Taylor's inequality on
the convex domain gives
\begin{equation}\label{eq:eleven-taylor}
 f\le g+u(t-t_0)+v(r-r_0)
       -\tfrac32\bigl[(t-t_0)^2+(r-r_0)^2\bigr].
\end{equation}
Write
\[
 a=1-qt_0/d,\quad b=1+qt_0/d,\quad
 \lambda=b+r_0/M,\quad \mu=1+\beta t_0-r_0,\quad
 D=q/2+\beta,\quad C=M-Dt_0.
\]
Expanding $(1-x)\log(1-x)$ gives
\begin{align*}
 g(s)&=C\log s+C_0+\sum_{m\ge2}a_ms^{m-1},\\
 a_m&=\frac{M(a^m+b^m-\lambda^m)+2-\mu^m
       -\{(1+t_0)^m+(1-t_0)^m\}/2}{m(m-1)}<0.
\end{align*}
Indeed $a=1/2$, $\lambda>b$, $\mu>1$ and
$M/4+2<\mu^2+(1+t_0)^2/2$. The series converges throughout the required interval, and $g''\le-C/s^2$. Moreover,
\begin{align*}
 u'={}&-\frac Ds-\frac{1+t_0}{2[1-(1+t_0)s]}
       +\frac{1-t_0}{2[1-(1-t_0)s]}\\
 &+\frac{qb}{2(1-bs)}-\frac{qa}{2(1-as)}
       -\frac{q\lambda}{2(1-\lambda s)}-\frac{\beta\mu}{1-\mu s},\\
 v'={}&-\frac{\lambda}{1-\lambda s}+\frac{\mu}{1-\mu s}.
\end{align*}
Pairing the two $1\pm t_0$ terms, and the positive $b$ term with the negative $\lambda$ term, gives $u'<0$. Since $\lambda>\mu$, also $v'<0$ and $v''<0$. One more differentiation gives
\[
 u''\le D/s^2+\frac{(1-t_0)^2}{2[1-(1-t_0)s]^2}.
\]

\paragraph{Two size ranges.}
Here $C=3591/2000$, $D=601/200$, $\lambda=346/225$ and $\mu=2569/2000$.
For $s\le9/20$, the derivative bounds give $u''\le D/s^2+3/500$.
The differences
$-v'=\lambda/(1-\lambda s)-\mu/(1-\mu s)$ and
$-v''=\lambda^2/(1-\lambda s)^2-\mu^2/(1-\mu s)^2$
increase in $s$, so substitution at $s=9/20$ gives $|v'|<2$ and $|v''|<20$. Put
\[
 J(s)=g(s)+u(s)/10-3/200+v(s)^2/6.
\]
The endpoint evaluations use the following three derivative formulas. Write $A_0=(d-qt_0)/2$,
$B_0=(d+qt_0)/2$, and $y_0=\beta t_0-r_0$.
\begin{align*}
 u(s)={}&\tfrac12\log\frac{[1-(1+t_0)s](1-t_0)}
                              {[1-(1-t_0)s](1+t_0)}
   +\tfrac q2\log\frac{B_0(M-A_0s)}{A_0(M-B_0s)}\\
 &+\tfrac q2\log\frac{(A_0-r_0)[M-(B_0+r_0)s]}{s(r_0+qt_0)^2}
   +\beta\log\frac{(1-y_0)(1-\mu s)}{y_0^2s},\\
 v(s)={}&\log\frac{(A_0-r_0)[M-(B_0+r_0)s]y_0^2}
                       {r_0(r_0+qt_0)(1-y_0)(1-\mu s)},\\
 s^2g'(s)={}&\tfrac12\log\bigl((1-s)^2-t_0^2s^2\bigr)
   +M\log\frac{M[M-(B_0+r_0)s]}{(M-B_0s)(M-A_0s)}
   +\log(1-\mu s)-2\log(1-s).
\end{align*}
Substitution, with the logarithm bounds of Appendix~\ref{app:arithmetic},
gives the following rational enclosures. Use $J'=g'+u'/10+vv'/3$.
\[
\begin{array}{c|r@{\ <\ }c@{\ <\ }l}
 s& \multicolumn{3}{c}{\text{enclosures}}\\\hline
 87/200 & -.062 & J(s) & -.061\\
        & .022 & J'(s) & .024\\
 9/20   & .509 & u(s) & .510\\
        & -.027 & v(s) & -.026\\
 97/200 & -.079 & g(s) & -.078\\
        & -.038 & g'(s) & -.037\\
 1/2    & -.564 & u(s) & -.563\\
        & -.150 & v(s) & -.148
\end{array}
\]
For $s\le9/20$, monotonicity gives $u(s)>1/2$, so the first quadratic in~\eqref{eq:eleven-taylor} is maximized under $t\le1$ at $t=1$. Thus $f\le J(s)$. Also $v>-1/20$, and
\[
 J''\le-\frac{299}{200s^2}+\frac3{5000}+\frac53<0,
\]
since $(v'^2+vv'')/3\le5/3$. Its tangent at $87/200$ gives $J<-1/25$ on this interval.
For $9/20\le s\le1/2$, the tangent to $g$ at $97/200$ gives $g<-3/40$, while monotonicity gives $|u|<57/100$ and $|v|<3/20$. Hence
\[
 f\le g+(u^2+v^2)/6
 <-\frac3{40}+\frac{(57/100)^2+(3/20)^2}{6}
 <-\frac3{200}.
\]
Continuity includes the boundary of the parameter domain.
\end{proof}
\subsection{A perfect matching on a marked six-regular factor}\label{app:degree-seven}

\begin{lemma}\label{lem:degree-seven-marked}
In the marked six-regular pairing $P$, half the vertices have outdegree one and half have outdegree five. Let $M_0$ be an independent uniform perfect matching, and let $\sigma$ be positive on the outdegree-five class. Then a.a.s.
\[
 |\cut_{P\cup M_0}(S)|>\frac{21}{5}|\sigma(S)|
 \qquad(\varnothing\ne S\ne V(P)).
\]
\end{lemma}

\begin{proof}
By complementation and reversal it suffices, up to a factor two in the expected count, to consider a set of order $sn\le n/2$ containing $(1+t)sn/2$ vertices of outdegree one and $(1-t)sn/2$ of outdegree five, with $0\le t\le1$. Set $A=3-2t$, $B=3+2t$ and $\delta=1/5$. The set has $Asn$ outgoing and $Bsn$ incoming points. If $J$ counts the outgoing boundary arcs in $P$ and $K$ the matching boundary edges, then
\[
 |\cut_{P\cup M_0}(S)|=4tsn+2J+K.
\]
The deficient sets thus have $2J+K\le\delta tsn$.

\paragraph{Small sets.}
For a set of order $k\le C_0\log n$, deficiency gives $2J+K\le k/5$. Thus at least $9k/10$ of its outgoing points have partners inside the set, and the perfect matching has at least $2k/5$ internal pairs. There are at most $2^{6k}$ choices of outgoing points. A specified set of $a$ of them all have partners inside with probability at most $(C k/n)^a$ for an absolute $C$. Similarly, there are at most $2^k k^b$ ways to choose $b$ disjoint internal pairs and the chance of their occurrence is $\prod_{j=0}^{b-1}(n-2j-1)^{-1}$. By independence and a union bound the failure probability over these sets is at most
\[
 2\sum_{k\le C_0\log n}\bigl[C'(k/n)^{3/10}\bigr]^k=o(1)
\]
for any fixed $C_0$.

\paragraph{The joint count.}
For a fixed set the point pairing gives
\[
 \Pr(J=j)=\frac{\binom{Bsn}{Asn-j}\binom{(3-Bs)n}{j}}{\binom{3n}{Asn}},
 \qquad
 \Pr(K=k)=\frac{\binom{sn}{k}\binom{(1-s)n}{k}k!
                  (sn-k-1)!!\,((1-s)n-k-1)!!}{(n-1)!!},
\]
where the second expression is zero unless the two unmatched orders are even. These counts and the choices from the two vertex classes give the following exponent on the boundary $2J+K=\delta tsn$. Put $r=J/(sn)$ and $y=K/(sn)=\delta t-2r$, so $0\le r\le t/10$. With $L,H$ as above, the exponent is
\begin{align}
 F(s,t,r)={}&\tfrac12\{H(s(1+t))+H(s(1-t))\}\notag\\
 &+L(Bs)+L(3-Bs)+L(As)+L(3-As)\notag\\
 &-L((A-r)s)-L(rs)-L((4t+r)s)-L(3-(B+r)s)-L(3)\notag\\
 &+L(s)+L(1-s)-L(ys)-\tfrac12L((1-y)s)-\tfrac12L(1-(1+y)s).
 \label{eq:d7-entropy}
\end{align}
The two boundary rates are concave in their separate boundary sizes, with maxima at their limiting means. Twice the first mean plus the second, divided by $sn$, is $2A(1-Bs/3)+1-s\ge A^2/3+1/2\ge5/6>\delta t$. Moving a deficient pair toward the pair of means increases the exponent until it reaches the equality boundary. Stirling's formula contributes a uniform $O(\log n)$ error to the exponent $nF(s,t,r)$.

It remains to show $F(s,t,r)<-s/500$ on this boundary. The number of class-size and boundary patterns is polynomial in $n$, so this gives failure probability $o(1)$ for orders at least $C_0\log n$ with $C_0$ sufficiently large.

\paragraph{Curvature of the exponent.}
Put $f=F/s$. For a column vector $z$ write $z^2=zz^T$. The Hessian with respect to $(t,r)$ is
\[
 \nabla^2 f=(T-P)(1,0)^2-a_1(2,1)^2-r^{-1}(0,1)^2
            -(4t+r)^{-1}(4,1)^2-b_1(\delta,-2)^2,
\]
where
\begin{align*}
 T&=4\left(\frac1A+\frac1B+\frac{s}{3-As}+\frac{s}{3-Bs}\right),\\
 P&=\frac12\left(\frac1{1+t}+\frac1{1-t}
                  +\frac{s}{1-(1+t)s}+\frac{s}{1-(1-t)s}\right),\\
 a_1&=\frac1{A-r}+\frac{s}{3-(B+r)s},\\
 b_1&=\frac12\left(\frac1{1-y}+\frac2y+\frac{s}{1-(1+y)s}\right).
\end{align*}
The identity
\[
 T-P=h(t)+\frac{s}{1-s}h\left(\frac{st}{1-s}\right),
 \qquad h(z)=\frac{8/3}{1-4z^2/9}-\frac1{1-z^2}
\]
gives $T-P<17/5$. Indeed, $h(z)-17/10$ has the sign of
$-1/30+7z^2/30-34z^4/45$ on $[0,1)$, and this quadratic in $z^2$ has negative discriminant. Since $b_1\ge 1/y$, $r+y/4\le t/10$ and $4t+r\le41t/10$, Cauchy--Schwarz gives, for real $u,v$,
\begin{align*}
 \frac{v^2}{r}+b_1(\delta u-2v)^2+\frac{(4u+v)^2}{4t+r}
 &\ge \frac1t\left(10(2v-u/10)^2+\frac{10}{41}(4u+v)^2\right)\\
 &=\frac1t\left(\frac{39}{10}u^2+30v^2+\frac{21}{205}(u-10v)^2\right).
\end{align*}
Since $t\le1$, we can retain $25v^2$ and use $5v^2$ together with $a_1(2u+v)^2$ to obtain an additional $4u^2/(16/5-2t)$, since $a_1\ge(3-2t)^{-1}$. Applying Cauchy--Schwarz once more,
\[
 \frac{39}{10t}+\frac{4}{16/5-2t}
 \ge\frac{(\sqrt{39/5}+2)^2}{16/5}
 >\frac{(279/100+2)^2}{16/5}>\frac{17}{5}+\frac{15}{4}.
\]
Thus $-\nabla^2 f\succeq\operatorname{diag}(15/4,25)$ throughout the interior.

\paragraph{Scalar bounds.}
At a reference point $(t_0,r_0)$ write $g(s)=f(s,t_0,r_0)$, $u(s)=f_t(s,t_0,r_0)$, and $v(s)=f_r(s,t_0,r_0)$. The domain in $(t,r)$ is convex. The curvature inequality gives
\begin{equation}\label{eq:d7-tangent}
 f\le g+u(t-t_0)+v(r-r_0)-\frac{15}{8}(t-t_0)^2-\frac{25}{2}(r-r_0)^2 .
\end{equation}
In the following scalar calculations all parameters are evaluated at the reference point. Write $a=A/3$, $b=B/3$, $\lambda=(B+r_0)/3$, $\mu=1+y$, and $C=5/2-(21/10)t_0$. Expanding $(1-z)\log(1-z)$ in~\eqref{eq:d7-entropy} gives
\[
 g(s)=C\log s+C_0+
 \sum_{j\ge2}\frac{3(a^j+b^j-\lambda^j)+1-
    \{\mu^j+(1+t_0)^j+(1-t_0)^j\}/2}{j(j-1)}s^{j-1}.
\]
The reference points used below satisfy $t_0>3/4$, $0<r_0<t_0/10$, $a<1/2$ and $1<\mu<\lambda<2$. The numerator is negative since $3a^j+1<7/4$ and $(\mu^j+(1+t_0)^j)/2>65/32$. Thus $g''\le-C/s^2$. For the other derivatives, direct differentiation gives
\begin{align*}
 -u'(s)&=\frac{21}{10s}
 +\frac{t_0}{(1-(1+t_0)s)(1-(1-t_0)s)}
 +\frac{2a}{1-as}
 +\frac{2(\lambda-b)}{(1-bs)(1-\lambda s)}
 +\frac{\delta\mu}{2(1-\mu s)},\\
 -v'(s)&=\frac{\lambda}{1-\lambda s}-\frac{\mu}{1-\mu s}.
\end{align*}
In particular $u,v$ decrease, $v''<0$, and $u''\le21/(10s^2)$.

For $s\le1/4$ use $(t_0,r_0)=(.970,.018)$ and $\alpha=1-t_0=.030$. At $s=1/4$ the bounds
\[
 .5269<u<.5271,\quad .0870<v<.0872,\quad
 -.03614<g<-.03613,\quad .6988<g'<.6990
\]
follow from the logarithm bounds in Appendix~\ref{app:arithmetic}. The derivative formulas give $-11.299<u'<-11.298$ and $-1.187<v'<-1.186$ there. For all $s\le1/4$ we then have $u>(15/4)\alpha$, $v>0$, and $|v'|<6/5$. Maximizing~\eqref{eq:d7-tangent} subject to $t\le1$ gives $f\le J=g+\alpha u-15\alpha^2/8+v^2/50$. Now $J''\le-2/(5s^2)+(6/5)^2/25<0$. The displayed endpoint bounds give $J(1/4)<-.021$ and $J'(1/4)>.35$, so $J(s)<-.021$ on this interval.

For $s\ge1/4$ unrestricted maximization in~\eqref{eq:d7-tangent} gives $f\le g+2u^2/15+v^2/50$. The table bounds $g\le G$, $|u|\le U$ and $|v|\le V$. Use the tangent to $g$ at the interval midpoint, and use the endpoint values of $u,v$. All terminating decimals are exact.
\[
\begin{array}{c|cc|rrr}
[s_-,s_+]&t_0&r_0&G&U&V\\\hline
{[.250,.295]}&.980&.0190&-.0134&.285&.062\\
{[.295,.330]}&.960&.0180&-.0090&.224&.033\\
{[.330,.365]}&.938&.0175&-.0080&.181&.035\\
{[.365,.400]}&.905&.0170&-.0081&.181&.052\\
{[.400,.435]}&.865&.0160&-.0084&.186&.056\\
{[.435,.470]}&.820&.0150&-.0089&.205&.083\\
{[.470,.500]}&.775&.0145&-.0089&.203&.054
\end{array}
\]
Every row gives $G+2U^2/15+V^2/50<-1/500$. Continuity includes the boundary of the domain.

\paragraph{How to check the table.}
Use~\eqref{eq:d7-entropy} for $g$. The following grouped derivative formulas reduce the number of logarithms in the checks:
\begin{align*}
 g'&=-\frac{1}{2s^2}\log\frac{(1-as)^6(1-bs)^6(1-s)^2}
 {(1-\lambda s)^6(1-\mu s)(1-(1+t_0)s)(1-(1-t_0)s)},\\
 u&=\frac12\log\frac{(1-t_0)(1-(1+t_0)s)}{(1+t_0)(1-(1-t_0)s)}\\
 &\quad+2\log\frac{B(3-As)(A-r_0)[3-(B+r_0)s]}{(3-Bs)A(4t_0+r_0)^2s}
\\
 &\quad+\frac1{10}\log\frac{(1-y)(1-(1+y)s)}{y^2s},\\
 v&=\log\frac{(A-r_0)y^2[3-(B+r_0)s]}{r_0(4t_0+r_0)(1-y)(1-(1+y)s)}.
\end{align*}
For example, the sixth row has midpoint $181/400$. The logarithm bounds of Appendix~\ref{app:arithmetic} give $g(181/400)<-.0090822$ and $|g'(181/400)|<.006289$. Hence the tangent is less than $-.0090822+.0175(.006289)<-.0089$. At the two endpoints, the bounds for $u$ are $(.18944,.18946)$ and $(-.20418,-.20416)$, and those for $v$ are $(.08230,.08232)$ and $(-.01897,-.01894)$. These give $U=.205$ and $V=.083$. The same substitutions verify the other rows.
\end{proof}

\section{The degree-five gap}\label{app:p1strict}

\begin{proposition}\label{prop:p1repair}
Put $\delta_5=1/12002$. Then a.a.s. $G_{n,5}$ is connected and has a balanced signing $\sigma$ with $q_G(\sigma)\ge3+\delta_5$.
\end{proposition}

The proof has four steps.
\begin{enumerate}
\item We bound the sizes of the monochromatic components in a marked Jaeger orientation.
\item We prove expansion of the bipartite graph obtained by contracting these components.
\item We recolour vertices in large trees and unicyclic components, using a reserve of singleton trees to preserve balance. Every final monochromatic component is a tree of order at most $C_1$.
\item We verify weighted Hall inequalities for the final component graph. Large sets inherit expansion from step~2, and small sets are controlled by crossing pairs left unexposed during recolouring. Proposition~\ref{prop:p1matching} then gives the cut-ratio bound.
\end{enumerate}
The constants are as follows.
\begin{center}\small
\begin{tabular}{@{}lll@{}}\toprule
Symbol&Value&Role\\\midrule
$r$&$.099n\le r\le .101n$&number of trees on each shore before recolouring\\
$b$&$1029/1000$&exponential moment: $\sum_Tb^{|T|}\le .14n$ on each shore\\
$C$&$2000$&trees of order greater than $C$ are freed\\
$\ell$&$10^{-20}$&$\lfloor\ell n\rfloor$ singleton trees of each type are reserves\\
$f$&$\le3\ell n$&number of free vertices\\
$D_0$&$3C+2=6002$&degree bound for fixed trees\\
$C_1$&$6C+2=12002$&order bound for final trees\\
$D$&$3C_1+2=36008$&degree bound for final trees\\
$\delta_5$&$1/C_1$&cut-ratio margin used in the Hall inequalities\\
$\alpha$&$10^{-14}$&small sets after recolouring have at most $\alpha n$ trees\\\bottomrule
\end{tabular}
\end{center}
We use the following terms. The trees and unicyclic components of the marked orientation are \emph{old}. A vertex is \emph{free} if it lies in an old tree of order greater than $C$, in an old unicyclic component, or in a reserve; the other old trees are \emph{fixed} and keep their colours. A crossing edge is \emph{exposed} if it meets a free vertex, and \emph{hidden} otherwise. A fixed tree is \emph{touched} if it meets an exposed edge, and \emph{untouched} otherwise. A reserve is \emph{usable} if it has five distinct fixed neighbours, each meeting just one exposed edge.

\subsection{The marked model and its components}

Let $X$ count the Jaeger orientations of a uniform pairing of five labelled points per vertex. The orientations direct each pair of points, so a loop has two point orientations. Choosing uniformly among pairing-orientation pairs gives the marked law, with density $X/\E X$ on pairings. Half the vertices have type $A$ with outdegree four, and half have type $B$ with outdegree one. Fix the types and the exceptional point at each vertex, incoming in $A$ and outgoing in $B$. A uniform bijection matches the outgoing points to the incoming points.

Call an arc pure if its endpoints have the same type. Pure indegree in $A$ and pure outdegree in $B$ are at most one, so each pure component is a tree or has one directed cycle. Tree edges point away from the root or cycle in $A$ and towards it in $B$. Contract these components to form the crossing multigraph $H$. A tree of order $s$ has degree $3s+2$ and a unicyclic component has degree $3s$.

Reveal the pure arcs. The two crossing pools are now independent uniform bijections. Only tree roots have an unmatched exceptional point, so the $B$-to-$A$ pool is a perfect matching between the two shores of tree roots, with $r$ pairs. The $A$-to-$B$ pool has $L=3n/2+r$ pairs. The hypergeometric variable $r$ counts matches between two prescribed sets of $n/2$ points in the original bijection of size $5n/2$. Its mean is $n/10$ and variance $O(n)$, so $.099n\le r\le .101n$ a.a.s.

\subsection{The tree profile}

Let $S$ be the total progeny of a Galton--Watson tree with offspring distribution $\operatorname{Bin}(4,1/5)$. On each shore of $H$, for any fixed $g$ with $|g(s)|=O(b^s)$, where $b=1029/1000$,
\begin{equation}\label{eq:d5-profile}
 \frac{1}{n}\sum_{T\text{ a tree on that shore}}g(|T|)
 \ \xrightarrow{\mathbb P}\ \frac{1}{10}\mathbb E g(S).
\end{equation}
To see this, explore pure arcs forward from a vertex of type $A$ or backward from a vertex of type $B$, and let $Z_v$ be the number of reached vertices, including $v$. For a search stopped at $s\le R=\lceil200\log n\rceil$, each query finds a new pure descendant with conditional probability at most $.20001$ for large $n$. The numerator is at most $n/2$ and at most $4R$ pairs have been revealed. Reaching $s$ needs $s-1$ successes in at most $4(s-1)$ queries. The generating-function bound at $4/3$ gives
\[
 \Pr(Z_v\ge s)\le \left[\tfrac34(1+.20001/3)^4\right]^{s-1}
 \le .971^{\,s-1}.
\]
The same bound holds conditional on $v$ being a root, which specifies only the partner of the exceptional point. In particular $\max_v Z_v\le R$ a.a.s. The expected number of vertices on pure cycles is bounded by $2\sum_{j\ge1}(4/5)^j$, by counting the cyclically ordered vertices and their compatible points. Thus a.a.s. there are $O(\log^2 n)$ vertices in unicyclic components, since each is reached from a cycle vertex.

A vertex is a root with probability $1/5$. Conditional on being a root, a search with a fixed cutoff converges to the branching process, since each new partner has limiting success probability $1/5$ and revisits have probability $o(1)$. For a fixed cutoff $M$, switching two partners changes only $O(M)$ root summands. Besides changes of root status, a change requires an unchanged path of length less than $M$ to the first changed arc. Each changed tail has at most $M+1$ such ancestors, by the pure indegree bound in the search direction. Bounded differences give deviations $o(n)$ in probability. The normalized expected contribution from $M<Z_v\le R$ is $O(\sum_{s>M}(.971b)^s)$, tending to zero as $M$ increases. Since a tree root reaches its whole tree, this proves \eqref{eq:d5-profile}.

Write $F(z)=\mathbb E z^S$. The branching recursion gives $F(z)=z(4+F(z))^4/625$. The mean is $5$, and the variance is $80$, by the recursion for the first two moments, using offspring mean $4/5$ and variance $16/25$. Thus $\mathbb E S^2=105$. In particular the sum of squared crossing degrees in $H$ is at most $203n$ a.a.s., since $\mathbb E(3S+2)^2/5=201.8$ and the cyclic contribution is $o(n)$. Also there are at least $n/100$ singleton trees on either shore, since $\Pr(S=1)=(4/5)^4$.

For the size cutoff take $C=2000$. Iterating the branching recursion gives $F(b)\le 4/3$, since $b(16/15)^4<4/3$. We have
\[
 \log b\ge (b-1)/b>1/36,\qquad
 C\log b\ge\frac{58000}{1029}>56>24\log10 .
\]
Here $\log10<7/3$ follows from the first seven terms of the exponential series. For $s\ge C+1$ the sequences $s/b^s$ and $s^2/b^s$ decrease, since $(1+1/s)^2<b$. By \eqref{eq:d5-profile}, the total order and the sum of squared crossing degrees of the trees larger than $C$, together with all unicyclic components, are at most $10^{-21}n$ and $3\cdot10^{-17}n$, respectively, a.a.s. Indeed, their limiting normalized contributions are bounded by $(4/15)(C+1)10^{-24}$ and $(20/3)(C+1)^2 10^{-24}$, using $d_H(T)\le 5|T|$. We also have $\sum_T b^{|T|}\le .14n$ over the old trees of each shore, by \eqref{eq:d5-profile} and $F(b)/10\le 2/15$.

\subsection{Hall expansion before recolouring}

For an old tree $T$, write $D_T=3|T|+1$ for its degree in the larger pairing pool. Let $G(z)=zF(z^3)$ be the limiting generating function for $D_T$. The limiting second moment is $\E(3S+1)^2=976$. Hence by \eqref{eq:d5-profile}, on either shore $\sum_T D_T^2\le1000r$ with probability $1-o(1)$. Also $L/r\to16$ in probability, so we may take $15.9r\le L\le Mr$ with $M=16.01$. We show that, with probability $1-o(1)$, on either shore every set of $x\le .51r$ old trees has more than $1.025x$ opposite tree neighbours, and every set with $x\le .00102r$ has more than $1.2x$.

Condition on the old monochromatic components. Fix $X$ of $x$ old trees on one shore and $Y$ of $y=\lfloor\rho x\rfloor$ old trees on the other. If $Y$ contains all the tree neighbours of $X$, all root mates of $X$ lie in $Y$, with probability $\binom yx/\binom rx$. Independently, all large-pool points at $X$ avoid the $c$ points at opposite old trees outside $Y$. Write $A_X=\sum_{T\in X}D_T$. The latter probability is at most $(1-c/L)^{A_X}$, by sequential exposure without replacement. We allow matches into old cyclic components. Summing over $X,Y$ cancels $\binom rx$ from the root probability, leaving $\binom ry\binom yx$ times the average of the large-pool probability over independent uniform $X,Y$.

\paragraph{Small sets.}
Write $u=x/r$. For $u\le .005$, designate four large-pool points at each tree of $X$. The possible partners in $Y$ and in cyclic components number at most $\sqrt{1000ry}+O(\log^2 n)$ by Cauchy--Schwarz. Each of the $4x$ queries then has conditional success probability at most $2.1\sqrt{y/r}$, since $L-4x\ge15.8r$ and $y\ge1$. Put $h(v)=-v\log v-(1-v)\log(1-v)$. Using $\binom ry\le(er/x)^{\rho x}$ and $y h(x/y)\le\rho x h(1/\rho)$ gives the first-moment bound
\[
 \bigl[e^{\rho(1+h(1/\rho))}(2.1)^4\rho^2u^{2-\rho}\bigr]^x.
\]
The exponential factor is $(e\rho)^\rho(\rho-1)^{1-\rho}$. For $\rho=41/40$ it is less than $(451/160)(6/5)<4$, using $e<11/4$ and $40\cdot451/160<200<(6/5)^{40}$, the latter because $(6/5)^5>2$. For $\rho=6/5$ it is less than $(33/10)(9/5)<6$, since $5\cdot33/10<(9/5)^5$. The first bracket is therefore at most $4(2.1)^4(41/40)^2(.006)<1/2$ for $u\le.005$. For the second bracket and $u\le .00102$, the fifth power of the upper bound is $(6(2.1)^4(6/5)^2)^5(.00102)^4<(.7)^5$. Summing in either case, with a split at $\sqrt r$, gives $o(1)$.

\paragraph{Larger sets.}
Take $\rho=41/40$ and $.005\le u\le .51$. Put $z=1-c/L$. Since every old tree has at least four large-pool points,
\[
 0\le z\le1-4(1-\rho\cdot.51)/M<.882 .
\]
The recursion for $F$ implies $F(z^3)\le .62z^3$ on this interval. Indeed, iteration from zero stays below $.62z^3$ because $(.8+.2\cdot.62\cdot.882^3)^4<.62$. Consequently
\[
 \log G(z)
 =4\log z+4\log(.8+.2 F(z^3))
 \le4\log z+4\log .8+Bz,\qquad B=.483,
\]
using $\log(1+s)\le s$ and $.62(.882)^2<B$. For either shore the same bound, with an additional $o(1)$, holds for $\log G_{\rm emp}(z)$, where $G_{\rm emp}(z)=r^{-1}\sum_T z^{D_T}$. Indeed the tree-size frequencies converge at each fixed size by \eqref{eq:d5-profile}, and the limiting frequencies sum to one. The sum of absolute frequency differences is therefore $o(1)$ in probability, and it bounds the difference between $G_{\rm emp}(z)/z^4$ and $G(z)/z^4$ uniformly on $0<z\le1$. Since $G(z)/z^4\ge\Pr(S=1)>0$, taking logarithms preserves this uniform convergence.

We use the following consequence of Maclaurin's inequality. For a uniform sample of $k$ entries from $r$ nonnegative numbers, the expected product is at most the $k$th power of their arithmetic mean. This follows by smoothing two unequal numbers with fixed sum in the elementary symmetric sum. Thus averaging $z^{A_X}$ over $X$ gives at most $G_{\rm emp}(z)^x$. For any fixed $q\in(0,1)$, put $a=-M\log q$. Since $c=L(1-z)\le Mr(1-z)$, the preceding logarithmic bound gives
\[
 G_{\rm emp}(z)^x
 \le q^c\exp\left\{r\sup_{v>0}
 \left[a+4u\log v+4u\log .8-(a-Bu)v\right]+o(r)\right\}.
\]
The case $z=0$ has zero probability and is included by continuity. Averaging $q^c$ over $Y$ by the same product inequality, now on the complement of $Y$, contributes at most $\exp\{(r-y)\log G(q)+o(r)\}$. Both empirical errors are uniform on the stated range. If $a>Bu$, the supremum occurs at $v=4u/(a-Bu)$. The first-moment exponent divided by $r$ is therefore at most $R(u)+o(1)$, where
\[
 R(u)=h(\rho u)+\rho u h(1/\rho)+(1-\rho u)\log G(q)+a
 +4u\left[\log\frac{16u}{5(a-Bu)}-1\right].
\]
For fixed $q$, this is convex where $a>Bu$ and $u\le.51$. The second derivative of the first term is $-\rho/[u(1-\rho u)]$, whereas that of $4u\log[u/(a-Bu)]$ is at least $4/u$, and $\rho/(1-.51\rho)<4$.

Choose $p\in(0,1)$ rational and define $q^3=625p/(4+p)^4$. Then $F(q^3)=p$, because $s\mapsto625s/(4+s)^4$ is strictly increasing on $[0,1]$. Hence $\log G(q)=\log p-a/M$. In the following table each choice of $p$ applies throughout the indicated interval. The last two columns are rational upper bounds on the endpoint rates.
\[
\begin{array}{cc|c|cc|rr}
 u_0&u_1&p&a_-&a_+&R(u_0)\text{ bound}&R(u_1)\text{ bound}\\\hline
 .005&.03&.875&.172161&.172162&-.0058&-.0063\\
 .03&.10&.725&.508596&.508597&-.0096&-.0116\\
 .10&.23&.565&1.103899&1.103900&-.0055&-.0073\\
 .23&.38&.44&1.845666&1.845667&-.0054&-.0077\\
 .38&.51&.365&2.479303&2.479304&-.0033&-.0031
\end{array}
\]
For the hand checks, enclose $a=-(M/3)\log[625p/(4+p)^4]$ by the logarithm formula in Appendix~\ref{app:arithmetic}. With $s=1-\rho u$ and $v=(\rho-1)u$, bound the five terms
\[
 -s\log(s/p),\quad -u\log u,\quad -v\log v,\quad
 (1-s/M)a,\quad 4u\left[\log\frac{16u}{5(a-Bu)}-1\right].
\]
Use $a_+$ in the fourth term and $a_-$ in the last. For instance, at $u=.51$, $p=.365$, the second term is less than $.51(.693148-2/101)<.343407$. Here $\log(.51)=\log(51/50)-\log2$, and $2/101$ is the first term of the series for $\log(51/50)$. The other four upper bounds are $-.127970$, $.055619$, $2.405398$, and $-2.679595$. The sum is less than $-.0031$. The denominators $a_--Bu_1$ are positive. Convexity gives a negative rate on the whole range, and summing over $x$ proves the assertion on either shore.

\subsection{The exposed graph and the recolouring}

Return to the exposure in which only the pure arcs are known. We use the Hall expansion estimate after recolouring, without conditioning the crossing pools on it. Free all vertices in old trees of order greater than $C$ and in unicyclic components. Put $\ell=10^{-20}$ and choose $\lfloor\ell n\rfloor$ singleton trees of each type as free reserves, in a fixed label order. All other trees keep their colours. Put $D_0=3C+2=6002$, the upper degree bound for fixed trees. Write $f$ for the number of free vertices. The profile estimates give $f\le3\ell n$ and a sum of squared degrees less than $4\cdot10^{-17}n$ for the free component nodes, since each reserve contributes $25$.

Expose the crossing edges incident to free nodes of $H$, and write $J$ for the graph of these exposed edges on the old component nodes. In each pool, after revealing precisely these pairs, the remaining bijection is uniform. The two remaining bijections are conditionally independent.

\begin{lemma}\label{lem:d5-pseudoforest}
With probability $1-o(1)$, the graph $J$ is a pseudoforest, that is, each of its components has at most one cycle, counting two parallel edges as a cycle of length two. Moreover, a component of $J$ containing an old unicyclic node is a tree and contains only one such node.
\end{lemma}

\begin{proof}
Consider a prescribed pattern with $v$ vertices and $e=O(\log n)$ edges. Conditional on the component data, each specified compatible point pair has probability at most $100/n$ at each step, since both pools have size at least $n/20$. Sum over distinct images of its vertices and choices of points. The sum of the point-choice factors at a vertex of degree $j\le4$ is at most $\sum_T d_H(T)^j=O(n)$ by the tree profile. For degree two it is at most $203n$, or $4\cdot10^{-17}n$ if the vertex is free. No edge of $J$ has two fixed ends. In a pattern formed by a bounded number of paths between endpoints and branch vertices, at least $v/2-O(1)$ degree-two vertices are free. Summing over at most $2^v$ free-status assignments therefore gives
\[
 O\bigl(n^{v-e}\lambda^v\bigr),\qquad
 \lambda=200\sqrt{203\cdot4\cdot10^{-17}}<1/2.
\]
The bounded number of exceptional vertices affects only the implied constant. This rules out paths of length $\lceil K\log n\rceil$ for large fixed $K$. A minimal connected subgraph with two independent cycles is either three internally disjoint paths with common endpoints or two cycles joined at a vertex or by a path. Each has $e=v+1$, and only polynomially many choices of path lengths for each $v=O(\log n)$, so the sum is $O(1/n)$. For a path joining two old unicyclic nodes, the endpoint factors are $O(\log^2 n)$ instead of $O(n)$, giving $O(\log^4 n/n)$. For a cycle with a path to an old unicyclic node, possibly on the cycle, $e=v$ and the designated node contributes only a power of $\log n$ instead of $n$. The sum is again $o(1)$.
\end{proof}

Expand each free node of $J$ back to its old component, retaining its pure edges. Leave each fixed tree contracted to a node of prescribed colour, and denote the exposed graph by $J_1$. Lemma~\ref{lem:d5-pseudoforest} implies that $J_1$ is a pseudoforest: expanding a tree preserves cycle rank, and an old unicyclic node lies in a tree of $J$ containing no other such node. On simple outcomes the graph $J_1$ is loopless. We perform the construction on loopless $J_1$, without imposing simplicity on the hidden pairs.

\paragraph{Balancing reserves.} Each usable reserve centres an isolated star in $J_1$ and can take either colour. For this estimate, query its five partners and then the remaining partners at those neighbours. In at most $5+5D_0$ queries each free hit has probability at most $125f/n$, since each pool still has at least $n/25$ available points. Repetitions among the five neighbours cost at most $250D_0/n$. The failure probability is thus less than $10^8\ell+10^7/n<1/16$ for large $n$. A switch affects only reserves at the changed endpoints or adjacent to the affected fixed nodes, so the failure count is $O_C(1)$-Lipschitz. Bounded differences give at most $\ell n/8$ unusable reserves per shore a.a.s.

\paragraph{Recolouring.} Usable reserves keep their original colour for the moment. The other components of $J_1$ are coloured by the following rules; fixed nodes always keep their prescribed colours.
\begin{enumerate}
\item In a tree component, choose a root; if the root is free, give it an arbitrary colour. Give every free non-root vertex the colour opposite to that of its parent.
\item In a unicyclic component whose cycle contains a fixed node, start at that node and give each free cycle vertex the colour opposite to that of its predecessor. If the cycle is entirely free, alternate the colours around it; on an odd cycle this leaves exactly one edge whose ends have equal colours. Colour the trees attached to the cycle by rule~1, keeping the assigned colours of their cycle roots.
\end{enumerate}
We check that every monochromatic component is a tree of order at most $C_1=6C+2$. No two fixed nodes are adjacent in $J_1$. A fixed node can have the colour of a free neighbour only if that neighbour is its parent or its cycle predecessor, so it is joined within its monochromatic component to at most one free vertex. The only possible monochromatic edge between free vertices is the equal-colour edge on an all-free odd cycle. Thus the monochromatic components are trees. Such a component consists of a lone fixed node, a single free vertex together with at most five fixed nodes, or of an equal-colour pair on a cycle together with at most six fixed nodes, since each vertex of the pair has two incidences on the cycle. Expanding the fixed trees, each of order at most $C$, gives the order bound $C_1=6C+2$.

The temporary change in the size of either colour class is at most $10^{-21}n+\ell n/4<\ell n/2$. Flip the necessary number of usable reserves of the appropriate colour to restore balance. These flips retain the order bound. A usable reserve has at most five fixed neighbours in its monochromatic tree in either colour. Make all choices by fixed rules, so the resulting two-forest partition is determined by the exposed data. Every hidden edge still joins opposite colours, since its endpoints are fixed vertices.

\subsection{Hall's inequality after recolouring}

Contract the final monochromatic trees to obtain the bipartite multigraph $H'$. A tree of order $s$ has degree $3s+2$. The two shores have the same number $r'$ of nodes, since each contains $n/2$ original vertices and the degree sums agree. The maximum degree is at most $D=3C_1+2=36008$. The exposed edges form a pseudoforest $Q$, obtained by contracting the final trees in $J_1$. With probability $1-o(1)$, the number $c_Q$ of its cycles is at most $\log n$. Indeed the bound on patterns with as many edges as nodes gives a bounded expected number of crossing cycles of length $O(\log n)$, and there are no longer cycles on the no-long-path event. The expected number of old pure cycles is also bounded. Expanding the free nodes adds only these pure cycles, and contracting trees cannot increase the number of independent cycles. Markov's inequality gives the claim.

At most $6f$ old nodes are free or adjacent to a free node. Each final tree other than an untouched old tree contains a free vertex or a touched fixed tree, so there are at most $6f$ such trees. Thus $|r'-r|\le6f$ and $.098n\le r'\le .12n$. On either final shore the sum of squared orders is at most $11n$ with probability $1-o(1)$. The untouched trees contribute at most $(10.5+o(1))n$ by the tree profile, and the remaining trees contribute at most $6fC_1^2<n/4$.

\paragraph{Small sets.}
Put $\alpha=10^{-14}$ and designate five incident points at each final node, using only the exposed data. For a set $X$ of $1\le x\le\alpha n$ nodes on one shore, pad a neighbourhood of size at most $2x$ to a set $Y$ of size $2x$ on the other shore. Let $U=X\cup Y$, $s=3x$, $m=e_Q(U)$, and let $c$ count the unicyclic components of $Q[U]$. There are $s-m$ tree components. Choose one root from each, and for each unicyclic component use a root on its cycle fixed in advance in $Q$. Order the roots by their labels and choose the component sizes with sum $s$. Each component admits a rooted plane spanning tree. There are at most $4^{a-1}$ rooted plane trees of order $a$ and at most $D$ choices of an edge of $Q$ at each parent. Since the spanning trees have $m-c$ edges in total, this gives at most
\[
 \binom{2r'}{s-m}\binom{c_Q}{c}2^s(4D)^{m-c}
\]
choices of $U$ with these parameters. Once $U$ and the shore of $X$ are fixed, $X$ and $Y$ are determined.

At most $m$ designated points of $X$ are exposed in a witness, so at least $5x-m$ designated hidden points must match into $Y$. In the larger pool there are at most $2Dx$ possible partners, and in the root pool at most $12x$, since a final tree contains at most six fixed old trees and hidden points belong to fixed trees. Even after the exposure of at most $5f$ pairs and the successive queries of at most $5x$ hidden points, the respective denominators are at least $1.49n$ and $.08n$. Thus each required match costs a factor at most $50000x/n$. The selected points are all on the same final shore, so no two can be partners. This estimate is conditional on the exposed data alone.

Write $j=s-m$. The binomial bound and $j\log(s/j)\le s/e$ give $\binom{2r'}j\le e^{s/e}(.66n/s)^j$, also at $j=0$, since $e<11/4$ and $2r'\le .24n$. Now $5x-m=j+2s/3$ and $.66\cdot50000/3<4D$. Multiplying the count by the matching probability therefore gives
\[
 \binom{c_Q}c \bigl[3\cdot10^8(s/n)^{2/3}\bigr]^s.
\]
Here $2e^{1/e}<3$, $(50000/3)^{2/3}<660$, and $3\cdot4D\cdot660<3\cdot10^8$. For the exponential comparisons, summing $1/j!$ through $j=4$ and bounding the tail by $1/100$ gives $8/3<e<11/4$. Also $e^{3/8}\le 1+3/8+(3/8)^2/(2(1-1/8))<3/2$. For $s\le\sqrt n$, summing the cycle binomials gives at most $(1+\log n)^s$, and the combined bracket is $o(1)$. For $s>\sqrt n$ the cycle sum is at most $2^{\log n}\le n$, and the bracket is below $.3$, since $(3\cdot10^8)^3(3\alpha)^2<(.3)^3$. In both ranges, summing over $m\le s$ and then $s$ gives $o(1)$. Thus every such $X$ has at least $2x$ neighbours a.a.s., uniformly in the good exposed data.

\paragraph{Intermediate and large sets.}
For $\alpha n\le x\le r'/2$ let $X_0$ consist of the untouched trees in $X$, so $|X_0|\ge x-6f$. At most $6f$ of its old neighbours are touched, and each untouched neighbour survives in $H'$. When $x\le10^{-4}n$, we have $|X_0|/r\le .00102$, so the old $1.2$ expansion gives
\[
 |N_{H'}(X)|\ge 1.2x-13.2f\ge1.19x.
\]
For $x\ge10^{-4}n$, use $|X_0|\le r'/2\le .51r$ and the old $1.025$ expansion to get $|N_{H'}(X)|\ge1.025x-13f\ge1.024x$. These comparisons use $f\le3\cdot10^{-20}n$.

Let $Y=N_{H'}(X)$ and $y=|Y|$. For an intermediate set with $y<2x$, the total order of the $u$ untouched trees in $Y$ is at most
\[
 36u\log(.14n/u)\le36y\log(.14n/y)<1176y .
\]
For $u=0$ the left side is zero. The first bound follows from Jensen's inequality and the exponential moment $\sum_T b^{|T|}\le .14n$ on the old shore, with $b=1029/1000$ and $\log b>1/36$. The function $u\log(.14n/u)$ increases up to $.14n/e$, so the second inequality follows from $u\le y\le2\cdot10^{-4}n$. For the last use $y\ge\alpha n$ and $\log10<7/3$. The touched trees contribute at most $6fC_1<\alpha n$, giving total order at most $1177y$.

\paragraph{Weighted Hall.}
Set $\delta_5=1/C_1=1/12002$ and give each tree $T$ demand $w(T)=1-\delta_5|T|/2$. Write $w(X)=\sum_{T\in X}w(T)$ for a set of trees. The demands are at least $1/2$ and their totals on the two shores agree. If $x\le r'/2$ and $y\ge2x$, then $w(Y)\ge x\ge w(X)$, which covers the small sets. For the intermediate sets with $y<2x$ we have
\[
 w(Y)\ge(1-1177/(2C_1))y>\tfrac{19}{20}\cdot 1.19x>x.
\]
For $10^{-4}n\le x\le r'/2$ and $y<2x$, Cauchy--Schwarz bounds the total order of $Y$ by $\sqrt{11ny}\le\sqrt{22nx}$. Thus
\[
 w(Y)-w(X)\ge .024x-\tfrac{\delta_5}{2}\sqrt{22nx}>0,
\]
since $\sqrt{22\cdot10^4}<470$ and $470/(2C_1)<.02$.

Every nonempty set of at most half a shore has strictly more neighbours than nodes. A set of more than half a shore has more than $r'/2$ neighbours by applying the estimate to a subset of size $\lfloor r'/2\rfloor$. These two observations imply connectedness. In a component with at most $r'$ nodes, two shores of size at most $r'/2$ would each have to expand into the other, and one shore of size greater than $r'/2$ would force the other to exceed $r'/2$ as well. The graph before contraction is therefore connected. For $x>r'/2$ the opposite-shore complement of $Y$ has at most half a shore and all its neighbours lie outside $X$. Apply its weighted inequality and use equality of the total demands to get $w(Y)\ge w(X)$.

The old expansion event and the hidden-pair estimates can be intersected after the separate probability bounds. Proposition~\ref{prop:p1matching} now gives $q_G(\sigma)\ge3+\delta_5$. The resulting numerical flow gap satisfies
\[
 \frac{\delta_5}{2+\delta_5}=\frac1{2C_1+1}=\frac1{24005}>\frac1{24006}.
\]

\subsection{Transfer to the uniform graph}

The orientation-count estimates of Pra\l{}at and Wormald~\cite[Section~2]{PW} verify the hypotheses of small-subgraph conditioning for $X$, with cycle means $4^j/(2j)$ and corrections $(-1/5)^j$. The distributional form~\cite[Theorem~4.3]{W} gives $X/\E X\xrightarrow{\mathrm d}W_5$ with $W_5>0$ almost surely. Indeed the logarithm of the Poisson product converges absolutely almost surely because the corrections exceed $-1$ and $\sum_j4^j5^{-j}/(2j)<\infty$.

Let $\mathcal B$ be the event that the pairing gives a simple graph failing the conclusion of Proposition~\ref{prop:p1repair}. The proof gives $\Pr_{\rm marked}(\mathcal B)=o(1)$, that is, $\E[X\mathbf 1_{\mathcal B}]/\E X=o(1)$. Lemma~\ref{lem:marked-transfer}, applied to $X$ and $X\mathbf 1_{\mathcal B}$, gives $X>X\mathbf 1_{\mathcal B}$ a.a.s. under the uniform law, and hence $\Pr_{\rm uniform}(\mathcal B)=o(1)$. Simplicity has probability bounded away from zero, and conditioning on simplicity gives $G_{n,5}$. This proves Proposition~\ref{prop:p1repair}.

\section{A forest obstruction at degree five}\label{app:p1forest}

We prove Proposition~\ref{thm:p1forest} by weighting vertices in tree components of order at least four and summing over all partitions into two induced forests.

\begin{proof}[Proof of Proposition~\ref{thm:p1forest}]
Work in the $5$-regular pairing model. For an ordered partition
$\mathcal P=(A,B)$ into two induced forests, let $L(\mathcal P)$ be the
number of vertices in components of order at least four. Fix $0<y<1$
and put
\[
 Z_y=\sum_{\mathcal P}y^{L(\mathcal P)}.
\]
Every partition with $L(\mathcal P)\le n/40$ contributes at least
$y^{n/40}$ to $Z_y$. Thus it suffices to prove that
$y^{-n/40}\mathbb E Z_y$ tends to zero for one choice of $y$.

\paragraph{The component series.}
A labelled tree $F$ has weight
$\prod_{v\in V(F)}(5)_{\deg_F(v)}$, where $(5)_t=5!/(5-t)!$ for
$0\le t\le5$ and $(5)_t=0$ for $t>5$. This counts its realizations
using the five labelled points at each vertex. Let $T(z)$ be the
exponential generating function of these weighted trees. Its terms
of orders at most three are
\[
 T_3(z)=z+\frac{25}{2}z^2+250z^3.
\]
In a planted tree, the root uses one of its five points for the parent
edge, leaving four available for children. Its generating function
$P(z)$ therefore satisfies $P=5z(1+P)^4$. Rooting a tree at a vertex gives $zT'(z)=z(1+P)^5$. Integration, with $T(0)=0$, yields
\begin{equation}\label{eq:p1-tree-series}
 T(z)=\frac P5-\frac{3P^2}{10}
      =z(1+P)^5-\frac{P^2}{2}.
\end{equation}
To give weight $y$ precisely to vertices in components of order at
least four, use
\[
 B(z,y)=T_3(z)+T(zy)-T_3(zy).
\]
The first term leaves the small trees unchanged, and the difference
weights every vertex of a larger tree by $y$.

\paragraph{Counting the two forests.}
On a fixed class of $k$ vertices, the total weight of forests with $r$
components is
\[
 \frac{k!}{r!}[z^k]B(z,y)^r.
\]
Such a forest uses $k-r$ internal edges and leaves $3k+2r$ points for
crossing edges. Hence class sizes $k_A,k_B$ and component counts
$r_A,r_B$ can contribute only if
\[
 k_A+k_B=n,\qquad
 m=3k_A+2r_A=3k_B+2r_B.
\]
The remaining points can then be paired across the classes in $m!$
ways. Choose $x>0$ with $B=B(x,y)<\infty$. Nonnegative coefficients
give $[z^k]B(z,y)^r\le x^{-k}B^r$, so the contribution of these counts
to $\mathbb E Z_y$ is at most
\begin{equation}\label{eq:p1-weighted-count}
 \frac{n!}{r_A!r_B!}\,B^{r_A+r_B}x^{-n}
       \frac{m!}{(5n-1)!!}.
\end{equation}
The factor choosing the classes cancels their two factorials; the
denominator is the total number of pairings.

\paragraph{A common bound for all counts.}
Put $\gamma_A=r_A/n$, $\gamma_B=r_B/n$ and
$w=m/n=3/2+\gamma_A+\gamma_B$. Stirling's formula bounds
\eqref{eq:p1-weighted-count} by $\exp\{nf+O(\log n)\}$, uniformly in
the counts, where
\[
 f=\gamma_A\log\frac B{\gamma_A}
   +\gamma_B\log\frac B{\gamma_B}
   +w\log w-\log x-\frac52\log5.
\]
Use the continuous value at a zero coordinate. For a fixed sum
$\gamma_A+\gamma_B$, concavity of $z\log(B/z)$ makes the first two
terms largest when the two component counts agree. Thus, with
$\gamma=(\gamma_A+\gamma_B)/2$, every partition is bounded by
\[
 f_B(\gamma)=2\gamma\log\frac B\gamma
       +\left(\frac32+2\gamma\right)
          \log\left(\frac32+2\gamma\right)
       -\log x-\frac52\log5.
\]
For $0<B<1/2$, differentiation gives
\[
 f_B'(\gamma)=2\log\frac{B(3/2+2\gamma)}\gamma,
 \qquad
 f_B''(\gamma)=-\frac3{\gamma(3/2+2\gamma)}<0.
\]
Its maximum over $\gamma>0$ is therefore attained at
$\gamma=3B/[2(1-2B)]$. Enlarging to this domain only increases the
bound. Substituting the maximizer and summing over the polynomially
many choices of counts gives
\begin{equation}\label{eq:p1-weighted-bound}
 \mathbb E Z_y\le n^{O(1)}
 \exp\left\{n\left[
 \frac32\log\frac3{2(1-2B)}-\log x-\frac52\log5
 \right]\right\}.
\end{equation}

Take $x=1/20$ and $y=125/324$. On the power-series branch,
$z=P/[5(1+P)^4]$ increases with $P$ for $0\le P<1/3$.
At $z=xy=25/1296$ this gives $P=1/5$, and
\eqref{eq:p1-tree-series} gives $T(z)=7/250$. Hence $B(x,y)<23/200$.
Combining \eqref{eq:p1-weighted-bound} with the initial bound
$y^{-n/40}\mathbb E Z_y$, the expected number of partitions with
$L(\mathcal P)\le n/40$ is at most
\[
 n^{O(1)}\rho^n,
 \qquad
 \rho=R^{1/2}y^{-1/40},\qquad R=\frac{432000}{456533}<\frac{19}{20}.
\]
In particular, $R^{20}<e^{-1}<3/8<125/324=y$, where $e>8/3$ follows from the first four terms of its series. Thus $\rho^{40}=R^{20}/y<1$. Markov's inequality excludes these partitions a.a.s. in the pairing
model. Conditioning on simplicity, using \eqref{s1:eq:simplicity},
proves the assertion for $G_{n,5}$.
\end{proof}

\needspace{28\baselineskip}
\section{Rational bounds for the scalar evaluations}\label{app:arithmetic}

All logarithms in the finite endpoint evaluations can be bounded by the
same calculation. For a positive rational $x$, write $x=2^k z$ with
$1\le z<2$. Put $w=(z-1)/(z+1)$. Integrating the geometric series for
$2/(1-w^2)$ gives
\[
 S(w)\le\log z\le S(w)+\frac{2w^{17}}{17(1-w^2)},
 \qquad S(w)=2\sum_{j=0}^7\frac{w^{2j+1}}{2j+1}.
\]
The same formula at $w=1/3$ bounds $\log2$. The interval width for
either normalized logarithm is at most $1/(68\cdot3^{15})$.
For $\log x=k\log2+\log z$, use the upper bound for $\log2$ in the
lower bound when $k<0$, and conversely. In a linear combination of
logarithms use the lower or upper endpoint according to the sign of the
coefficient. Terms $L(0)$ are zero and terms involving $\log e$ use
$\log e=1$. For example, the bounds at $w=1/3$ and $w=1/5$ give
\[
 .693147<\log2<.693148,\qquad
 1.098612<\log3<1.098613,
\]
which give $(793/500)\log2-\log3>7/10000$ in the bisection count.
All terminating decimals in the bounds are rational numbers.
The same substitution applies to the endpoint estimates in
Appendices~\ref{s3:app:small}, \ref{s3:app:near}, \ref{s3:app:tables}
and \ref{app:marked-expansion}.

All comparisons involving $\pi$ use only the classical bound $\pi<22/7$.

\section*{AI disclosure}
ChatGPT, including Codex, assisted with literature searches, proof
exploration, exposition, and arithmetic checks. An AI-generated referee
report also informed the revision. The proofs and the bounds for the
finite arithmetic checks are given in the paper.

\par\medskip\noindent\textit{Email address.}
\href{mailto:a.jafari@newuu.edu}{\nolinkurl{a.jafari@newuu.edu}}

\end{document}